\documentclass[a4paper,fleqn]{cas-sc}
\usepackage{graphics,graphicx,pdfpages}
\usepackage{CJK}
\usepackage{indentfirst}
\usepackage{bbm}
\usepackage{bm}
\usepackage{amsmath}
\usepackage{amssymb}
\usepackage{pdfpages}
\usepackage{float}
\usepackage{varioref}
\usepackage{float}
\restylefloat{table}
\restylefloat{figure}
\usepackage{caption}
\usepackage{subcaption}
\usepackage{fancyhdr}
\usepackage{algorithm}
\usepackage[noend]{algpseudocode}
\usepackage[section]{placeins}
\usepackage{stfloats}
\usepackage{natbib}

\UseRawInputEncoding
\def\tsc#1{\csdef{#1}{\textsc{\lowercase{#1}}\xspace}}
\tsc{WGM}
\tsc{QE}

\usepackage{xcolor}

\newcommand{\subjectto}{\text{\rm subject to}} 

\newcommand\norm[1]{\left\lVert#1\right\rVert}

\newcommand{\minimize}{\ensuremath{\mathop{\mathrm{minimize}}\limits}}

\begin{document}
\let\WriteBookmarks\relax
\def\floatpagepagefraction{1}
\def\textpagefraction{.001}

\shorttitle{Dynamic OD Estimation During COVID}    

\shortauthors{Yu et al.}  

\title [mode = title]{Extending Dynamic Origin-Destination Estimation to Understand Traffic Patterns During COVID-19}  
\begin{abstract}
Estimating dynamic Origin-Destination (OD) traffic flow is crucial for understanding traffic patterns and the traffic network. While dynamic origin-destination estimation (DODE) has been studied for decades as a useful tool for estimating traffic flow, few existing models have considered its potential in evaluating the influence of policy on travel activity. This paper proposes a data-driven approach to estimate OD traffic flow using sensor data on highways and local roads. We extend prior DODE models to improve accuracy and realism in order to estimate how policies affect OD traffic flow in large urban networks. We applied our approach to a case study in Los Angeles County, where we developed a traffic network, estimated OD traffic flow between health districts during COVID-19, and analyzed the relationship between OD traffic flow and demographic characteristics such as income. Our findings demonstrate that the proposed approach provides valuable insights into traffic flow patterns and their underlying demographic factors for a large-scale traffic network. Specifically, our approach allows for evaluating the impact of policy changes on travel activity. The approach has practical applications for transportation planning and traffic management, enabling a better understanding of traffic flow patterns and the impact of policy changes on travel activity.
\end{abstract}


\begin{highlights}
\item We extend the dynamic origin-destination estimation (DODE) formulation to capture traffic patterns in a large-scale urban setting.
\item We use our DODE method to analyze traffic sensor data from Los Angeles to understand travel patterns during the COVID-19 stay-at-home order in 2020.
\item We find that traffic flow generally decreased during the stay-at-home order, but it decreased more when the origin or destination was a higher income district.
\end{highlights}


\tnotemark[1] 

\tnotetext[1]{Research reported in this publication was supported by National Library of Medicine of the National Institutes of Health under award number 1R21LM013697. The content is solely the responsibility of the authors and does not necessarily represent the official views of the National Institutes of Health.}
%

\author[1]{Han Yu}[orcid = 0000-0001-5019-805X]



\ead{hyu376@usc.edu}



\author[1]{Suyanpeng Zhang}
\ead{suyanpen@usc.edu}

\author[1]{Sze-chuan Suen}
\ead{ssuen@usc.edu}

\author[1]{Maged Dessouky}
\ead{maged@usc.edu}

\author[2]{Fernando Ordonez}
\ead{fordon@uchile.cl}


\affiliation[1]{organization={Daniel J Epstein Department of Industrial and Systems Engineering},
            addressline={University of Southern California}, 
            city={Los Angeles},
            postcode={90089}, 
            state={CA},
            country={USA}}

\affiliation[2]{organization={Department of Industrial Engineering},
            addressline={Universidad de Chile}, 
            city={Santiago},
            country={Chile}}


\maketitle
\section{Introduction}
The COVID-19 pandemic has had a significant impact on transportation behavior, leading to widespread changes in how we move around our cities and communities. With restrictions on public transportation and an increase in remote work, traffic patterns have undergone significant changes, as shown in \cite{Tirachini2020COVID-19Needs}. Understanding these trends may require origin-destination traffic flow estimation \citep{Zhou2007AFramework}, which can provide important insights during this period of change. However, we find that traffic patterns obtained from classic OD estimation models fail to model correctly the day-long commuting traffic patterns expected during the pandemic. In this article, we explore the impact of COVID-19 on traffic and transportation using an extended OD model that accounts for some of these factors. These insights can then uncover how travel behavior can guide health policy decisions and vice versa.

The traffic pattern in Los Angeles (LA) County is influenced by a variety of factors, including population density, land use patterns, and transportation infrastructure. To better understand the traffic flow between different health districts in the county, it is necessary to estimate the origin-destination (OD) demand. However, accurately estimating OD demand for an extended period can be a challenging task. Although several studies have tackled the problem of dynamic OD estimation (DODE), obtaining time-dependent OD demand remains a complex and arduous process. Moreover, most research in this area has focused on estimations between two directions, for example, North to South \citep{Ma2018EstimatingData}. As a result, there is a significant gap in the literature when estimating dynamic OD demand for extended periods between different districts that cover a large urban area. OD demand may also vary across districts. There can also be significant differences between districts when estimating OD demand over an extended period. For instance, areas with high commercial activity may experience a surge in traffic during weekdays in commercial districts compared to residential areas. Most prior works develop a DODE model that formulates the temporal and spatial relationships between regions by creating a network of possible paths between regions (where links represent paths and nodes represent locations), but they often do not consider local versus highway traffic patterns or commuting behaviors. This omission can result in inaccuracies in predicting the distribution of traffic demand. To overcome this limitation, our approach explores network configuration effects and considers commuting and road size patterns through regularization terms to ensure more nuanced estimations.

Constructing an accurate network within the region of interest is another critical factor that can impact the accuracy and computing time of DODE. In DODE, networks are typically constructed to follow the road network, with the assumption that vehicles travel along the links in the network. Both the construction of the network and optimization model play a critical role in the quality of the OD estimates. A well-designed network can help reduce computational complexity and improve the estimation of OD demand. To address these issues, this paper proposes a data-driven approach that uses high-granular traffic flow counts and traffic speed data to estimate time-dependent OD demand. The proposed method leverages the DODE to predict OD demand over an extended period accurately. The approach takes into account the variation in OD demand between different districts and employs a well-designed network to optimize the computational complexity of the estimation process. Overall, the proposed method represents an advancement of OD demand estimation and has the potential to provide valuable insights into traffic patterns of large urban networks.

We use our DODE model to analyze the transportation behavior in LA County using empirical data in 2019 and 2020 to understand the impacts of the stay-at-home order on traffic patterns across social dimensions (household income). In LA County, various measures were implemented to reduce the transmission of COVID-19. The initial phase of a stay-at-home order was put in place from March 19, 2020, until April 28, 2020. During this period, indoor gatherings and outdoor events with at least 10 people in a confined space were prohibited. Additionally, all indoor malls, shopping centers, playgrounds, and non-essential businesses were required to close. After April 28, 2020, restaurants and bars were permitted to reopen. However, due to a significant increase in cases, on July 1, 2020, the state ordered the closure of bars and indoor dining until August 28, 2020. On November 30th, as COVID-19 cases surged across the nation, Los Angeles County introduced a new stay-at-home order, known as the "safer at home" order. This order imposed additional restrictions on businesses and gatherings without completely banning them. It took effect until December 20. The order prohibited all public and private gatherings with individuals outside of a person's household. Residents were advised to stay home as much as possible. We observe a clear drop in traffic during the pandemic period, likely due to the implementation of social distancing measures and stay-at-home policies. Our DODE model provides more nuanced insights as it can identify which origin-destination pairs were most affected by the stay-at-home order and use these outcomes to understand social or economic trends in traffic behavior due to the pandemic guidelines. This information can be used by transportation planners to understand how public health rules differently impact different areas in terms of traffic. In particular, we examine the relationship between traffic demand and annual household income across districts in LA. Our analysis reveals that traffic distribution is affected by income, where higher income areas responded more to stay-at-home orders than lower income areas. Overall, our proposed DODE optimization model provides valuable insights into traffic behavior over heterogeneous regions. 

The remaining sections are organized as follows. In Section \ref{sec:model}, we discuss the formulation of our approach. In Section \ref{sec:numericalresult}, we propose a framework for selecting the appropriate network and parameters and present the results of a real-world experiment where we estimate OD flow from 2019 to 2020 in LA County, and explore the relationship between annual income and OD flow. Finally, we draw conclusions in Section \ref{sec:conclusions}.

\subsection{Related Literature}
Our research has  enhanced DODE estimation methods, enabling a more precise measurement of the impact of COVID-19 policies on people's travel patterns, particularly for the road network. In this section, we will delve into relevant literature from three key perspectives: DODE estimation, road network assessment, and the influence of Covid-19 on traffic patterns.
\subsubsection{DODE estimation}
We first begin by providing context around DODE estimation, then provide background on prior work studying traffic patterns during COVID-19. DODE estimation has been an area of study for several decades, with most prior work relying on either a general least squares model or a state-based model.  \cite{Hendrickson1984ESTIMATIONREGRESSION.} used regression models to estimate congestion data. Then \cite{Bell1991TheSquares} and \cite{Cascetta1993DynamicCounts} extended the basic regression model to generalized least square (GLS) models to estimate dynamic OD demands. Within several decades, GLS-based models have been widely used to estimate traffic flow demand, including single-level and bi-level GLS models, with the latter considering the effect of congestion. A bi-level optimization approach integrates the conventional GLS estimation model and the standard network equilibrium model into one process. At first,  \cite{Yang1992EstimationNetworks} used a bi-level model, where the upper-level model minimizes deviation, and the lower level performs dynamic traffic assignment for static OD estimation. Then \cite{Zhou2003DynamicApplications} extended the bi-level optimization model to incorporate multi-day traffic data. Regarding its application to congested networks, \cite{Frederix2013DynamicPractice} explored bi-level optimization with the goal of alleviating congestion. Meanwhile, \cite{Shao2014EstimationCounts} proposed a bi-level optimization model aimed at estimating traffic demand during peak hours. While most of the work focused on static estimation, \cite{Cantelmo2014AnDemand} presented an in-depth analysis of the bi-level gradient approximation approach for dynamic demand. However, the general bi-level formulation for OD estimation is non-continuous and non-convex, which limits its scalability. To address this, \cite{Lundgren2008AProblem} proposed a heuristic method to solve the bi-level OD estimation problem. On the other hand, a single-level formulation is proposed for OD estimation incorporating stochastic user equilibrium constraints. \cite{Nie2008ADemands} and \cite{Nie2010ATables} formulated a single-level DODE model with user equilibrium solved by a variational inequality. Recently, \cite{Ma2018StatisticalData} developed a generalized single-level formulation for OD estimation that incorporates stochastic user equilibrium constraints. The DODE model also has benefited from numerous machine learning techniques, such as those proposed by \cite{Cao2021Day-to-dayData} and \cite{Fan2022DeepData}, to enhance its capabilities.

To estimate OD, traffic counts are one of the most commonly used sources of data, which are usually collected at strategic locations in a transportation network. \cite{Zhang2008EstimatingAnalysis} evaluated the roles of count data, speed data, and history OD data in the effectiveness of DODE. Recently, GPS and probe vehicle data have also been used to estimate traffic flow. \cite{Munizaga2014ValidatingData} used GPS data to track the movement of individual vehicles in a transportation network, while probe vehicle data, which is collected from GPS-enabled devices in vehicles, can be used to estimate traffic flows between origins and destinations based on travel time and speed measurements as introduced in \cite{Antoniou2006DynamicSources}. Although there are many prior studies examining DODE, none consider the origin and destination traffic demand based on district level data, and scalability of the estimation remains a challenge, particularly when applied to real-world traffic data. In this paper, we examine methods to do this for analysis of COVID-19 transportation patterns.
\subsubsection{Covid-19 Policy and traffic}
We are interested in pandemic traffic patterns as COVID-19 mitigation policies have had a profound impact on commuting patterns. Previous research conducted to analyze the COVID-19 pandemic has largely focused on analyzing raw flow or volume data, rather than examining the impact of the pandemic on OD patterns. Measures such as stay-at-home orders, social distancing, and remote work have led to changes in travel behavior and traffic patterns \citep{Abdullah2020ExploringPreferences}. In particular, the pandemic has affected OD traffic flow, as the typical movement of people and goods between regular locations has been disrupted. Several studies have examined the impact of COVID-19 on OD traffic flow in specific regions or cities. \cite{Sevtsuk2021TheMA} examined the impact of COVID-19 on amenity visits in American cities using aggregate app-based GPS positioning data from smartphone users, which typically constitute over half of all urban trips. \cite{Harantova2020ComparisonCOVID-19} examined the flow and speed of vehicles in the cities Zilina and Cadca, Slovakia, during the pandemic period, and found COVID-19 reduced the flow of traffic and had a positive effect on air pollution. During the COVID-19 pandemic, public transit was significantly impacted as well. \cite{Brough2021UnderstandingPandemic} and \cite{Mashrur2022WillChoice} documented the magnitudes of and mechanisms behind socioeconomic differences in travel behavior during the COVID-19 pandemic, and discovered that the COVID-19 pandemic caused a significant decrease in travel and use of public transit. \cite{Orro2020ImpactSpain} focused on the main origin-destination flows of urban mobility based on the city bus network in Spain and reported a decrease in bus ridership. Their findings suggest that this reduction is not uniform across all bus networks. \cite{Bian2021TimeSeattle} studied the time lag in COVID-19 policy making. They proposed a Bayesian change point detection method to quantify the time lag effect reflected in transportation systems when authorities took action in response to COVID-19.

However, there is a limited number of studies that analyze the flow of vehicle traffic between different origins and destinations in order to understand mobility patterns. \cite{Tao2023RevealingChina} analyzed vehicle traffic flow from origins to destinations in Nantong, China, based on taxi trajectory data. This analysis has limitations as taxi usage represents only a small portion of overall travel behavior and may not be the primary mode of transportation for most people.  \cite{Gu2021AnalysisCOVID-19} used data on cellular networks to estimate the traffic flow between two groups of city pairs around Hubei Province in China during COVID-19. While the focus of their study was on intra-city travel patterns, our study focuses on investigating inter-city travel patterns,  taking into account shorter distances and perhaps more complex routes and road networks. Despite the considerable research on COVID-19 and traffic pattern changes, to the best of our knowledge, we are not aware of any prior work that considered the OD traffic flow on vehicles in urban areas to study the impact of COVID-19 on commuting patterns.

\subsection{Contributions}
This manuscript extends the DODE literature by addressing challenges that naturally arise when estimating origin-destination outcomes on real-world traffic data across districts and measure the impact of Covid-19 policy on traffic pattern. First, we modify the DODE formulation of \cite{Ma2018EstimatingData} by incorporating additional regularization terms to account for daily traffic patterns (e.g., traffic flow that exits and returns to a district due to commuting patterns) and interplay between local and arterial road traffic. These regularization terms may result in outcomes that better reflect observed commuting behavior. Secondly, because real-world traffic data typically operates on a complex road network, particularly in urban areas, we examined the impact of road network modeling choices on the DODE outcomes. This analysis demonstrates that the outcomes may be highly sensitive to network configuration and parameters, and we provide guidance on how appropriate selections can be made. Finally, to demonstrate this methodology, we apply our framework to estimate DODE outcomes in Los Angeles County (LAC) between 2019 and 2020, when the COVID-19 pandemic resulted in drastic changes to OD patterns. This numerical experiment on a large-scale network using real-world data shows that our framework can help extend DODE estimation to complex urban settings to generate policy insights. To this end, we additionally examine our DODE outcomes by demographic subgroup to understand how stay-at-home policies impacted different populations. We integrate income data from the census to examine traffic flow changes between 2019 and 2020 (likely due to stay-at-home orders) for DODE pairs with high income origin or destination districts. Our results suggest that OD pairs involving at least one high income district were more likely to have decreased traffic flow in 2020 during the height of the COVID-19 pandemic.

\section{The Model}\label{sec:model}
\begin{figure}
        \centering
\includegraphics[width=0.9\linewidth]{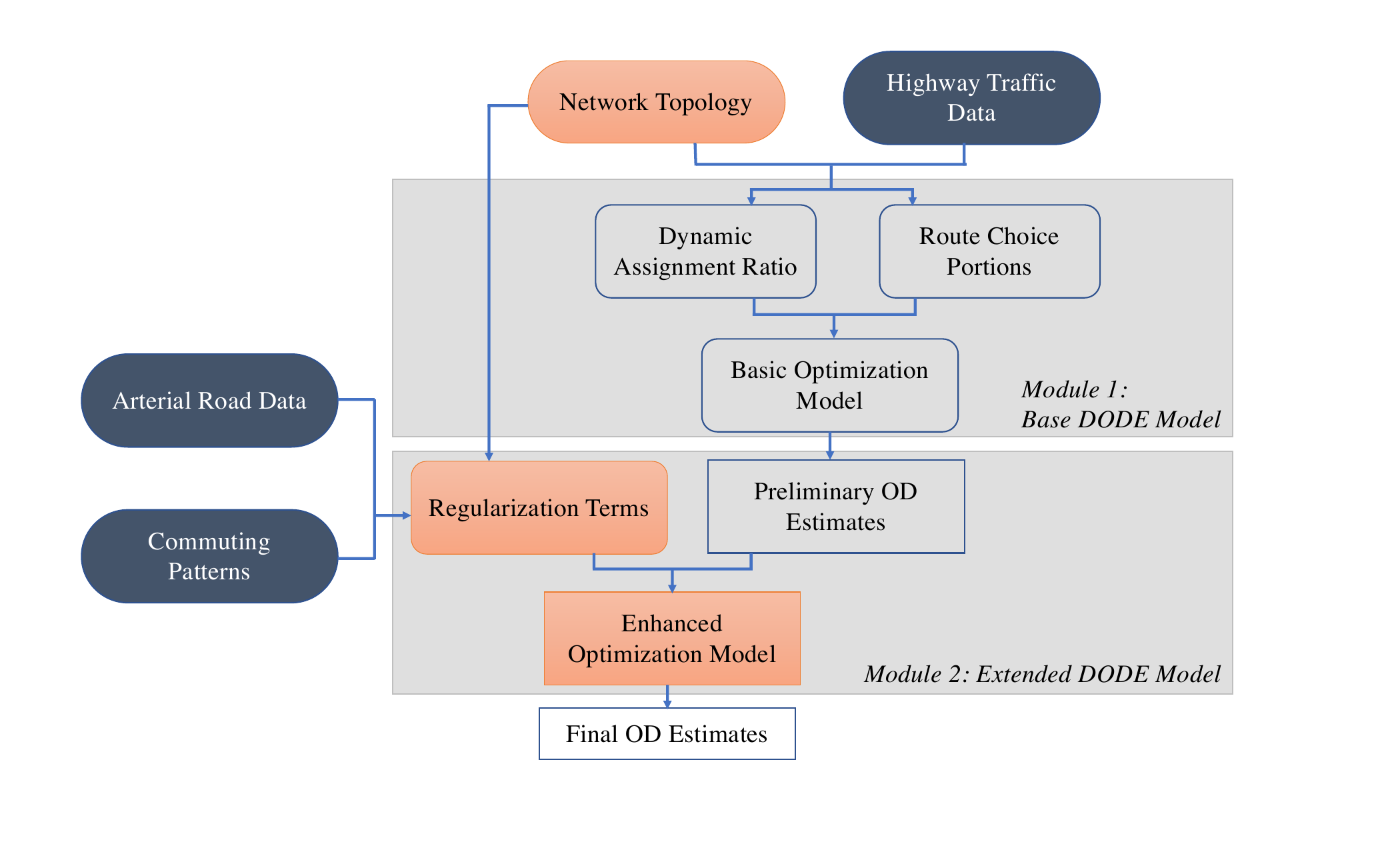}
	\caption{General Framework.}\label{fig:general_framework}
\end{figure}

In this section, we introduce a model to utilize the high-granular traffic demand and speed data to estimate 24/7 dynamic OD traffic flow. The objective of our study is to estimate traffic flow between regions by leveraging continuous traffic volume and speed data gathered through sensor monitoring. The proposed method comprises two major modules, as depicted in Figure \ref{fig:general_framework}. The first module, the base DODE model, employed by the \citet{Ma2018EstimatingData} model to obtain preliminary OD estimates. The second module, extended DODE model, refines the estimates by incorporating constraints to improve accuracy and realism in the estimation of traffic flow between health districts.

The first module is based on prior work from \citet{Ma2018EstimatingData}, as described in Section \ref{sec:basic_model}. \citet{Ma2018EstimatingData} presented a model incorporating dynamic assignment ratio (DAR) to represent the evolution of traffic flow in discrete time, as well as a route choice model to determine path selection. They formulated the dynamic origin-destination estimation as a non-negative least squares problem. These models use a traffic network, where each node represents a possible origin, destination, or intermediate points, and the shortest path between each origin-destination pair is determined exclusively using roads. We therefore establish a network based on the physical road network, which enables us to monitor traffic flow using sensor data. Our model utilizes this data to estimate DODE, as described in more detail in section \ref{sec:network_intro}.

In the second module, we extend the \citet{Ma2018EstimatingData} model by adding several additional regularization terms to capture traffic behavior in a large-scale urban context. To distinguish pass-through flow and origin-destination flow, we leverage local road data. Since local roads feed the highway network, local traffic volume should be loosely proportional to origin and destination volume. We therefore incorporate a lower bound on origin and destination volume as a function of local road traffic. This ensures that origin and destination nodes located downtown with complex local arterial roads have a minimum traffic threshold for the total amount of traffic flowing in and out of the node. 

In this work, we are interested in estimating daily traffic patterns. Workplace commuting patterns mean that people tend to travel back to their origin region by the end of the day. To incorporate this knowledge into our model, we have implemented a symmetry constraint to ensure that the inflow and outflow of traffic from each region are similar.

These modifications may result in large deviations in total volume for each node from the base OD model. We therefore also include a regularization term to restrict our outcomes to have similar node volumes to the original model, as these should be unchanged due to commuting patterns or local road information.

In the next sections, we describe each module, and its components, in detail.

\subsection{Network Definition}\label{sec:network_intro}

\begin{figure}[H]
        \centering
\includegraphics[width=0.7\linewidth]{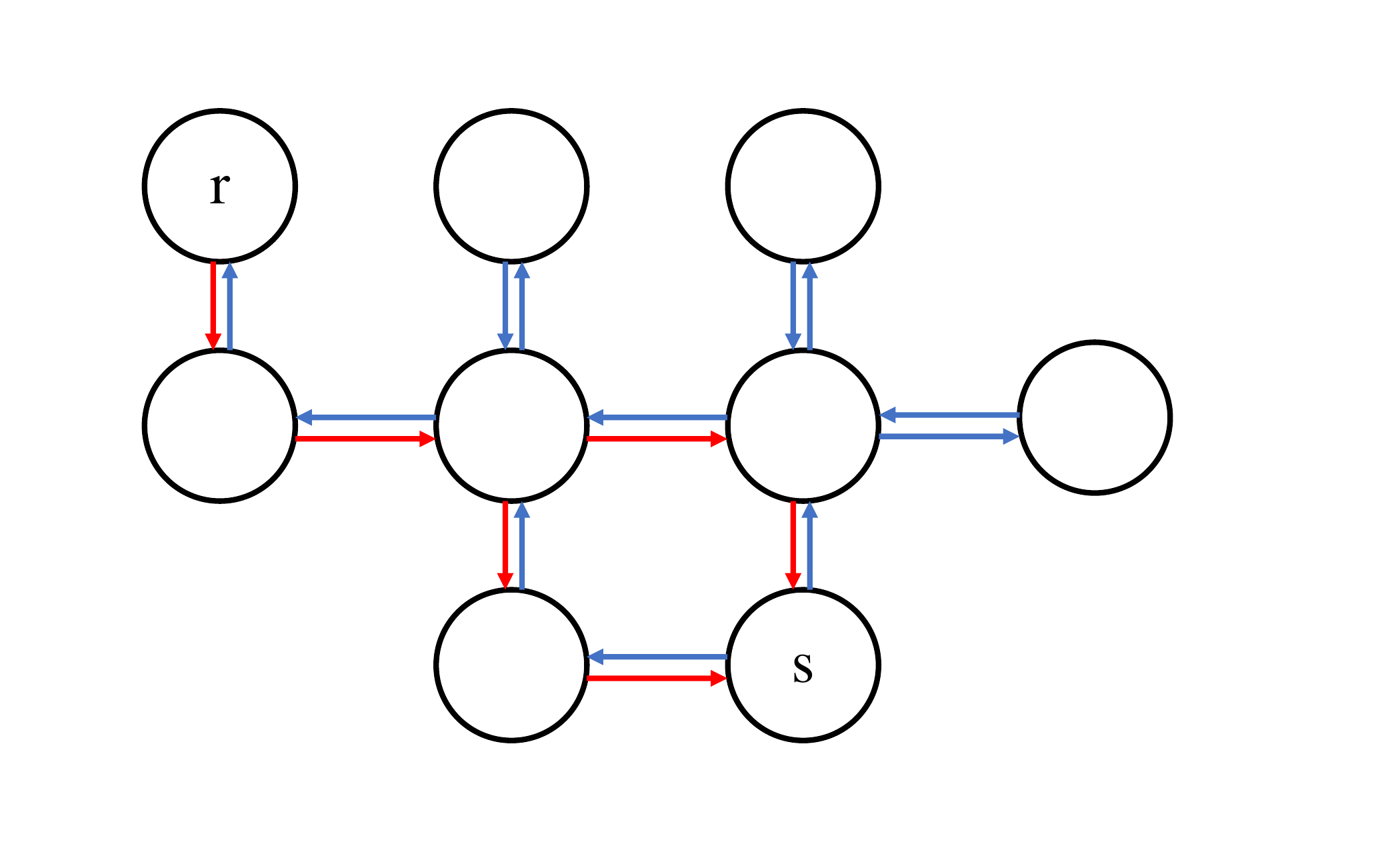}
	\caption{Example graph $G$ with OD $rs$ and links.}
        \label{fig:sample_network}
\end{figure}

The DODE model relies on a traffic network, which we describe here. A traffic network can be thought of  as a graph that represents the physical road network, where the nodes of the graph represent intersections or other points of interest, and the edges represent the roads that connect them. Let us consider a road network, denoted by $G$, that consists of a set of nodes $N$, a set of links $L$, and a set of OD pairs, denoted by $\Psi$. An OD pair denoted by ${r,s}$ with origin $r$ and destination $s$ belonging to node set $N$, is connected by links. We assume a connected graph, that is, for every OD pair $r$, $s$, there exists a path from the origin node $r$ to the destination node $s$. Figure \ref{fig:sample_network} shows an example origin Node $r$ and destination node $s$, with links representing roads between regions shown in blue. Red arrows indicate the shortest paths from $r$ to $s$. In our analysis, if there are multiple shortest paths, one is selected at random for consideration in the DODE model. In the DODE model, our objective is to estimate the traffic flow between each OD pair, which is assumed to flow through the selected path. We obtain the traffic flow data for each link from real-world observations, then utilize the DODE model to estimate the traffic flow for each OD pair. 

Network structure plays an important role in traffic flow analysis, as the design of the traffic network can greatly impact the predicted traffic flow, while the size of the network can affect the computation speed. In our case, we aim to investigate the traffic flow between different origin and destination pairs, with a focus on people traveling from one location to another by vehicle. Therefore, the network should follow the road network, and the number of nodes that denote the origin and destination should be appropriately selected. In section \ref{sec:ne_network}, we conducted a numerical experiment to select an appropriate traffic network in LA County.


\subsection{General Framework}

The aim of our study is to estimate day-to-day origin-destination (OD) traffic flows using continuous traffic flow and speed records collected by sensors. In this section, we present the general framework of the proposed method.

Our objective is to minimize the difference between the observed link flow and the estimated path flow. We define the OD traffic flow as $q$ and the link flow as $y$, and assume that there is a relationship between them, denoted by $\hat{y} = f(q)$, where $\hat{y}$ is the estimated link flow. The function $f(\cdot)$ maps the OD flow to the link flow. Our objective function can therefore be formulated as equation \eqref{eq:loss}, where we aim to minimize the gap between the estimated link flow and the observed true link flow. We assume that there is a linear relationship between $y$ and $f$, which can be expressed as $f(q) = Aq$, where $A$ is a known matrix, which we describe below.

\begin{equation}
\label{eq:loss}
    \begin{array}{rl}
         &  \minimize_{\bm q\geq 0} \norm{\bm y-\hat{\bm y}}^2 \\
       =  & \minimize_{\bm q\geq 0} \norm{\bm y-f(\bm q)}^2 \\ 
       =  & \minimize_{\bm q\geq 0} \norm{\bm y-\bm A\bm q}^2 \\
    \end{array}
\end{equation}

 Our proposed method consists of two modules. The first module is a DODE model proposed by \cite{Ma2018EstimatingData}, which establishes the basic relationship between links and ODs. The second module is an extended model with several regularization terms to improve the estimation accuracy. The $A$ matrix therefore comprises these two parts: $A_{b}$ that captures the link and OD relationship, and $A_{e}$ to capture regularization terms. Thus define $A$ as $A = [A_{b}; A_{e}]$. If we only considered the first module from \cite{Ma2018EstimatingData}, equation \eqref{eq:loss} would be simply formulated as $\min_{q\geq 0} \epsilon_{b}^2$, where $\epsilon_{b} = \norm{y-A_{b}q}$. A detailed description of $A_{b}$ is provided in Section \ref{sec:basic_model}.

In the second model, to improve the estimation of OD flow obtained from the first model, we introduce several regularization terms into the objective using the same structure. We define the decision variable $x$ to include $q$ variables and slack variables. We incorporate local road information as a lower bound, denoted as $\epsilon_{LB}=\norm{b_{LB} - A_{LB} x}$, daily symmetric traffic pattern information as $\epsilon_{s}=\norm{b_{s} - A_{s} x}$, and total traffic flow as 
$\epsilon_{\tau}=\norm{b_{\tau} - A_{\tau} x}$ into the model. The objective can now be formulated as equation \eqref{eq:2ndformulation}. When the regularization parameters $\eta$, $\beta$, and $\gamma$ are set to zero, equations \eqref{eq:loss} and \eqref{eq:2ndformulation} are equivalent. A more detailed explanation of the second model is provided in Section \ref{sec:2nd_model}.
\begin{equation}
\label{eq:2ndformulation}
    \begin{array}{rl}
         &  \minimize_{q\geq 0} \norm{y-\hat{y}}^2 \\
        = & \minimize_{q,x\geq 0} \epsilon_{b}^2 + \eta\epsilon_{LB}^2 + \beta\epsilon_{s}^2 + \gamma\epsilon_{\tau}^2\\
       = & \minimize_{q,x\geq 0} \left\{\norm{y-A_{b}q}^2 + \eta\norm{b_{LB} - A_{LB} x}^2 + \beta\norm{b_{s} - A_{s} x}^2 + \gamma\norm{b_{\tau} - A_{\tau} x}^2\right\}
    \end{array}
\end{equation}

\subsection{Module 1: Base DODE Model}\label{sec:basic_model}

The first model utilizes a dynamic assignment ratio and route choice portions to estimate preliminary OD traffic flow from link flow. The dynamic assignment ratio accounts for the dynamic changes in traffic demand across time, while the route choice component considers the various possible paths that drivers can take to reach their destinations. The model used traffic volume and speed data to estimate one-day dynamic OD. We assume traffic average speed remains constant within each time interval and estimates the average traffic volume in each time interval across all links in our network. For all links $a\in L$, we can find a link flow $y_a$. The objective function of the problem computes the L2 norm error between the observed link flow $y_a(t)$ at time $t$ and the estimated link flow $\hat{y}_a(t)$. $T$ is the set of all time intervals. Then $\hat{q}^k_{rs}(t')$ represents the estimated traffic flow from origin $r$ to destination $s$ following path $k$ at time $t'$. Then the model can be formulated as equation \eqref{eq:twostage_RandD}.
\begin{equation}
\begin{array}{cl}
    \displaystyle \minimize_{\{\hat{q}^k_{rs}(\cdot )\}_{k,r,s}}
     &\quad\displaystyle \sum_a\sum_{t\in T}\norm{y_a(t)-\hat{y}_a(t)}^2\\
    \subjectto &\quad \hat{q}^k_{rs}(t') \geq 0\quad \forall t'\in T\;\;\forall \{r,s\}\in \Psi\;\;\forall k \in \Phi_{rs},
\end{array}
\label{eq:twostage_RandD}
\end{equation}
where the estimated link flow $\hat{y}_a(t)$ and estimated path flow $\hat{q}^k_{rs}(t')$ have relationship:
\[\hat{y}_a(t) = \displaystyle\sum_{rs\in \Psi}\sum_{k\in \Phi_{rs}}\sum_{t'\in T}\rho_{rs}^{ka}(t,t')\hat{q}_{rs}^k(t')\]

Here, $\Psi$ is the set of OD pairs, and $\Phi_{rs}$ is the set of paths from $r$ to $s$. $\rho_{rs}^{ka}(t,t')$ is the dynamic assignment ratio (DAR), which is the portion of the $k$th path flow departing within time interval $t'$ between OD pair $rs$ which arrives at link $a$ within time interval $t$. As we assumed that the vehicles are spread evenly in time and on each link of the network, $\rho_{rs}^{ka}(t,t')$ is calculated by the speed and volume data from each link, as shown in \cite{Ma2018EstimatingData}.

We next computed the route choice portions for all OD pairs. Route choice portions represent the probability of a driver choosing a particular path. We define the route choice portion $p_{rs}^k(t)$ as the probability that a vehicle chooses path $k$ between OD pair $rs$ at time $t$. It distributes the OD demand $q_{rs}(t)$ to the $k$th path flow $q^k_{rs}(t)$ using equation \eqref{eq:routechoice}. If we consider only the shortest path between each origin $r$ and destination $s$ pair, the route choice portion can be represented by an identity matrix.
\begin{equation}
    \label{eq:routechoice}
    p^k_{rs}(t) = \frac{q^k_{rs}(t)}{\displaystyle\sum_{k\in\Phi_{rs}}q^k_{rs}(t)}
\end{equation}
Then the final basic optimization model is:
\begin{equation}
\begin{array}{cl}
    \displaystyle \minimize_{\{\hat{q}^k_{rs}(\cdot )\}_{k,r,s}}
     &\quad\displaystyle\sum_{a\in L}\sum_{t,t'\in T}\norm{y_a(t)-\sum_{rs\in \Psi}\sum_{k\in \Phi_{rs}}\rho_{rs}^{ka}(t,t')p^k_{rs}(t')\hat{q}_{rs}^k(t')}^2\\
    \subjectto &\quad \hat{q}^k_{rs}(t') \geq 0\quad \forall t'\in T\;\;\forall r,s\in \Psi\;\;\forall k \in \Phi_{rs},
\end{array}
\label{eq:twostage_final}
\end{equation}
Equation \eqref{eq:twostage_final} can be reformulated in a more general form, as introduced in equation \eqref{eq:loss}. In this context, we define the decision variable OD traffic flow as $\bm q$, which is a vector of $q_{rs}^k(t)$ in the model, using the equations described in \eqref{eq:qy}. On the other hand, the input data, observed link flow $\bm y$, is defined in appendix equation \eqref{eq:qy} by element $y_a(t)$. The matrix $A_{b}$ can be formulated using elements $\rho_{rs}^{ka}(t,t')$ and $p^k_{rs}(t')$, as shown in Appendix \ref{ap:basic}. Through these transformations, we can reformulate equation \eqref{eq:twostage_final} as: 
\begin{equation}
    \minimize_{\bm q\geq 0} \norm{\bm y-\bm A_{b}\bm q}^2
    \label{eq:1stpart}
\end{equation}


\subsection{Module 2: Extended DODE Model}\label{sec:2nd_model}

In this section, we introduce three additional regularization terms that can be incorporated into the model to improve its accuracy and adaptability to more complex scenarios. The experiment studied by \cite{Ma2018EstimatingData} focused on a single pair of origin and destination with only two highways. However, in practice, traffic flow data often involve multiple origin-destination pairs and a larger road network with multiple highways. Therefore, we aim to extend the model to be applicable to more general and complex cases.

To accomplish this, we add three regularization terms that can account for variations in traffic flow across multiple OD pairs and highways. The first regularization term encourages a lower bound on the total OD flow of each origin and destination. This can help distinguish whether the flow originated or ended at a particular location and prevent the model from generating unrealistic or invalid estimates by ensuring pass-through flow is attributed correctly. The second regularization term constrains the inflow and outflow to be symmetric for each region. This follows our intuition that people tend to return to where they started within a day and can help to prevent the model from generating estimates that are inconsistent with this intuition. Finally, the third regularization term encourages the total traffic flow to be within a reasonable range. This can help to prevent the model from generating estimates that are excessively high or low. By incorporating these regularization terms into the model, we can improve its performance and extend its applicability to more complex real-world scenarios.

\subsubsection{Incorporating local road traffic information}
The base model does not have any constraints to distinguish pass-through flow and OD flow. For example, suppose location A connects to location B and location B connects to location C. It is difficult to distinguish whether there is only one vehicle traveling from location A passing location B to location C, or whether there are two vehicles one vehicle traveling from location A to location B and another vehicle traveling from location B to location C. To resolve this issue, we use the traffic flow information from local roads to provide a lower bound of total traffic flow at each node at time $t$. We then introduce positive slack variables $\bar{x}_i(t)$ to convert the inequality into an equation, as shown in Equation \eqref{eq:localroad_ineq}.
\begin{equation}
\begin{array}{rl}
    l_{i}(t) + d_{i}(t) &  \geq \alpha LB_i(t)\\
     l_{i}(t) + d_{i}(t) - \bar{x}_i(t)&  = \alpha LB_i(t)
\end{array}
\label{eq:localroad_ineq}
\end{equation}
We denote traffic flow $l_i(t)$ originating from node $i$ and traffic flow $d_i(t)$ that has a destination at node $i$. $LB_i(t)$ represents the lower bound on traffic flow that starts or ends at node $i$ at time $t$, which is calculated by local road traffic information, where $\alpha$ is the discount parameter. To obtain $LB_i(t)$ from real data, we add up the traffic volume measured by all sensors on arterial roads located within an $\lambda$-kilometer radius circle centered on node $i$, at time $t$, using the node's location (latitude and longitude) as a reference point. In the experimental section, we set $\lambda$ to be one kilometer. While this volume provides an estimate of the flow within this radius, it may overestimate the complete origin and end of traffic flow at node $i$. To address this limitation, we apply a discount parameter $\alpha$ to the lower bound $LB_i(t)$. Next, $l_r(t)$ and  $d_r(t)$ can be represented by $q$ as equation \eqref{eq:localroad_ld}. 
\begin{equation}
    l_{i}(t) = \sum_{\{is \in \Psi, s\in N\}}\sum_{k\in \Phi_{is}}q_{is}^k(t);\quad d_{i}(t) = \sum_{\{ri \in \Psi, r\in N\}}\sum_{k\in \Phi_{ri}}q_{ri}^k(t)
    \label{eq:localroad_ld}
\end{equation}
If we incorporate it as a regularization term in the first basic optimization model, equation \eqref{eq:localroad_ineq} can be expressed as a L2 norm loss, as shown in equation \eqref{eq:regularizer_LB}.
\begin{equation}
    \begin{array}{cl}
    \epsilon_{LB}^2&=\displaystyle \sum_{i\in N}\sum_{t\in T} \norm{l_{i}(t) + d_{i}(t) - \bar{x}_i(t) - \alpha LB_i(t)}^2\\
    & = \displaystyle \sum_{i\in N}\sum_{t\in T} \norm{\sum_{\{is \in \Psi, s\in N\}}\sum_{k\in \Phi_{is}}q_{is}^k(t) + \sum_{\{ri \in \Psi, r\in N\}}\sum_{k\in \Phi_{ri}}q_{ri}^k(t) - \bar{x}_i(t) - \alpha LB_i(t)}^2
    \end{array}
    \label{eq:regularizer_LB}
\end{equation}
In the next step, we can formulate equation \eqref{eq:regularizer_LB} as a general regularization term $\epsilon_{LB} = \norm{\bm b_{LB} - \bm A_{LB}\bm x}$ in equation \eqref{eq:2ndformulation}. The vector $\bm{x}$ consists of both the OD traffic flow $\bm{q}$ and slack variables $\bm{\bar{x}}$, while the vector $\bm{b}_{LB}$ contains elements $LB_i(t)$. The elements in $A_{LB}$ are either 1, -1, or 0, which correspond to the parameters in front of $\bm q$ and $\bar x$. These formulations can be expressed mathematically as shown in Appendix \ref{ap:local_road}.

By applying this transformation, we can rewrite equation \eqref{eq:regularizer_LB} as equation \eqref{eq:LB_final}, which allows us to incorporate the traffic flow constraints \eqref{eq:localroad_ineq} into the first base optimization model \eqref{eq:loss}. 
\begin{equation}
    \begin{array}{cl}
    \epsilon_{LB}^2&=\displaystyle \sum_{i\in N}\sum_{t\in T} \norm{l_{i}(t) + d_{i}(t) - \bar{x}_i(t) - \alpha LB_i(t)}^2\\
    & = \norm{\bm b_{LB} - \bm A_{LB}^q\bm q - \bm A_{LB}^{\bar{x}}\bar{\bm x}}^2\\
    & = \norm{\bm b_{LB} - \bm A_{LB}\bm x}^2
    \end{array}
    \label{eq:LB_final}
\end{equation}

This results in an enhanced model that can effectively optimize the traffic flow based on highways, while also satisfying the local road flow constraints. In summary, the regularization term provides a means of balancing the optimization between arterial roads and local roads, ensuring the total pass-through traffic flow is not overestimated and allowing the model to capture the complex interactions between these roads types, producing more realistic traffic flow solutions.

\subsubsection{Incorporating Daily Traffic Pattern Information: Symmetry Constraints}
In general, one would assume that the daily inflow and outflow of traffic from each region are roughly equal since people tend to travel to different regions during the day for work, school, or other activities and return to their starting point at night. Therefore, we expect that the input and output traffic flows across a day would be similar in magnitude.

Several nodes may be within a region. To define the regions in our model, we first created a set of regions $R = \{R_1, R_2, \ldots, R_n\}$, where each $R_i$ represents a distinct region. Each region may include several nodes. Each $R_i$ includes all the nodes within that district, so it is itself a set of $r_1, r_2, \ldots, r_m$, such that $\{r| r\in R_i\}$. We use $rs$ to define an OD pair that starts at node $r$ and ends at node $s$. In this way, we are able to organize the nodes in our model by region, which is important for capturing the geographical structure of the data. The total inflow (outflow) to a region is the sum of all inflows (outflows) across nodes within that region. 

We require that the total traffic flow $f_{R_i, R_j}$ from region $R_i$ to region $R_j$ should be similar to the traffic flow $f_{R_j, R_i}$ from region $R_j$ to region $R_i$ within a day. In other words, the amount of traffic flow between any two districts should be roughly the same in both directions:
\begin{equation}
     f_{R_i, R_j} = f_{R_j, R_i}
     \label{eq:sym_eq}
\end{equation}
where $f_{R_i, R_j}$ denotes the total traffic flow from region $R_i$ to region $R_j$, and $f_{R_j, R_i}$ denotes the total traffic flow from region $R_j$ to region $R_i$. By incorporating this constraint into our DODE model, we are able to capture more realistic commuting patterns between different regions. We define a set $D$ consisting of all the days in our dataset, denoted as $D = {d_1, d_2, \ldots, d_n}$. For each day $d_i \in D$, we further divide it into time intervals denoted by $t\in d_i$. In this way, we can represent the temporal structure of our data by organizing it into days and time intervals within each day. This organization allows us to capture the time-varying dynamics of the traffic flow and model how it changes over different days and times. We express the total traffic flow from region $R_i$ to region $R_j$, denoted as $f_{R_i, R_j}$, in terms of the $q_{rs}^k(t)$ variables. We can define $f_{R_i, R_j}$ using the following equation:
\begin{equation}
    f_{R_i, R_j} = \displaystyle\sum_{t\in d}\sum_{\{rs\in\Psi|r\in R_i, s\in R_j\}}\sum_{k\in \Phi^{rs}}q_{rs}^k(t)
\end{equation}

We can express our desire for symmetric inflows and outflows as:
\begin{equation}
    \begin{array}{cl}
    \epsilon_{s}^2&=\displaystyle \sum_{d\in D}\sum_{R_i,R_j\in R}\norm{f_{R_i, R_j} - f_{R_j, R_i}}^2\\
    & = \displaystyle \sum_{d\in D}\sum_{R_i,R_j\in R} \norm{\sum_{t\in d}\sum_{\{rs\in\Psi|r\in R_i, s\in R_j\}}\sum_{k\in \Phi^{rs}}\left[q_{rs}^k(t) - q_{sr}^k(t)\right]}^2
    \end{array}
    \label{eq:regularizer_sym}
\end{equation}

However, we wish to reformulate this into the form $\norm{\bm b_{s} - \bm A_{s}\bm x}$ for use in equation \eqref{eq:2ndformulation}. As in the previous section, the vector $\bm{x}$ includes both the OD traffic flow $\bm{q}$ and slack variables $\bm{\bar{x}}$. However, in this case, we do not need to use slack variables, so the corresponding parameters in the matrix $A_{s}$ are set to zero, and the vector $\bm{b}_{s}$ is also zero. The elements in $A_{s}$ take on values of either 1, -1, or 0, which correspond to the parameters in front of $\bm q$ and $\bar x$. By applying this transformation as shown in Appendix \ref{ap:sym}, we can rewrite equation \eqref{eq:regularizer_sym} as equation \eqref{eq:sym_final}, which allows us to incorporate the traffic flow symmetry constraints \eqref{eq:sym_eq} into the first basic optimization model \eqref{eq:1stpart}. 
\begin{equation}
    \begin{array}{cl}
    \epsilon_{s}^2&=\displaystyle \sum_{d\in D}\sum_{R_i,R_j\in R}\norm{f_{R_i, R_j} - f_{R_j, R_i}}^2\\
    & = \norm{\bm 0 - \bm A_{s}^q\bm q}^2\\
    & = \norm{\bm b_{s} - \bm A_{s}\bm x}^2
    \end{array}
    \label{eq:sym_final}
\end{equation}

\subsubsection{Imposing fidelity to total traffic flow}

The solution to the base DODE model (the outcome of module 1) does not incorporate these symmetry considerations for daily traffic patterns or consider local road information, but we would like to maintain some fidelity to the original solution so that our enhancements do not wholly determine the resultant outcome. The most direct way to achieve this is by constructing another error term which measures the deviation from the base solution and use it as a regularization term. By changing the regularization weights, we can directly control the relative importance of module 1 and module 2 outcomes. 

We define this error term as the sum of all traffic flow $q_{rs}(t)$ between the origin $r$ and destination $s$ for all time intervals $t$ within a day $d$. This value should be at a similar scale to the total flow in the first base optimization model, which does not include any additional constraints and is predicted solely by the observed highway link flow. We use $\hat{q}_{rs}^{k}(t)$ to denote the preliminary OD flow estimation in the first model. Here the constraint is defined as the following equation.
\begin{equation}
    \sum_{t\in d}\sum_{rs\in\Psi}\sum_{k\in\Phi^{rs}} q^k_{rs}(t) = \sum_{t\in d}\sum_{rs\in\Psi}\sum_{k\in\Phi^{rs}} \hat{q}_{rs}^{k}(t)
    \label{eq:total_eq}
\end{equation}
If we incorporate the constraint as a regularization term in the first basic optimization model, equation \eqref{eq:total_eq} can be expressed as an L2 norm loss, as shown in equation \eqref{eq:regularizer_total}.
\begin{equation}
    \begin{array}{cl}
    \epsilon_{\tau}^2&=\displaystyle \sum_{d\in D}\norm{\sum_{t\in d}\sum_{rs\in\Psi}\sum_{k\in\Phi^{rs}} \hat{q}_{rs}^{k}(t) - \sum_{t\in d}\sum_{rs\in\Psi}\sum_{k\in\Phi^{rs}}q^k_{rs}(t)}^2\\
    \end{array}
    \label{eq:regularizer_total}
\end{equation}
In the next step, we can formulate equation \eqref{eq:regularizer_total} as a general regularization term $\epsilon_{\tau} = \norm{\bm b_{\tau} - \bm A_{\tau}\bm x}$ in equation \eqref{eq:2ndformulation}. Similar to the previous section, the vector $\bm{x}$ includes both the OD traffic flow $\bm{q}$ and slack variables $\bm{\bar{x}}$ as defined in equation \eqref{eq:LB}. In this case, we do not need to use slack variables, so the corresponding parameters in the matrix $A_{\tau}$ are set to zero. The vector $\bm{b}_{\tau}$ contains elements $\hat{q}_{rs}^{k}(t)$. The elements in $A_{\tau}$ take on values of either 1, or 0, which correspond to the parameters in front of $\bm q$ and $\bar x$.

By applying this transformation, we can rewrite equation \eqref{eq:regularizer_total} as equation \eqref{eq:total_final} as shown in Appendix \ref{ap:total}, which allows us to incorporate the traffic flow symmetry constraints \eqref{eq:total_eq} into the first base optimization model \eqref{eq:1stpart}. 
\begin{equation}
    \arraycolsep=1.4pt\def\arraystretch{2.2}
    \begin{array}{cl}
    \epsilon_{\tau}^2&=\displaystyle \sum_{d\in D}\norm{\sum_{t\in d}\sum_{rs\in\Psi}\sum_{k\in\Phi^{rs}} \hat{q}_{rs}^{k}(t) - \sum_{t\in d}\sum_{rs\in\Psi}\sum_{k\in\Phi^{rs}}q^k_{rs}(t)}^2\\
    & = \norm{\bm b_{\tau} - \bm A_{\tau}^q\bm q}^2\\
    & = \norm{\bm b_{\tau} - \bm A_{\tau}\bm x}^2
    \end{array}
    \label{eq:total_final}
\end{equation}
When this regularization term is introduced into the first model, it ensures that the total traffic flow remains at a reasonable level. It can help prevent the model from becoming too sparse, meaning it avoids too many OD flow estimations becoming zero. This is important because a model that is too sparse may not capture all the relevant information in the data, leading to poor performance in the estimation of OD flows. By using a regularization term, the model is encouraged to strike a balance between the number of variables or features used and the accuracy of the estimations.



\section{Numerical Analysis: Los Angeles Traffic During the 2020 COVID Stay-at-home Order}\label{sec:numericalresult}
Using the enhanced framework described above to calculate DODE values allows us to better capture traffic behavior for a large-scale traffic network using real data. In this analysis, we examine traffic flows during the 2020 COVID-19 stay-at-home order. We demonstrate that using our DODE framework can provide more insight than examining traffic flow volumes alone, and also show how our solutions can be used to make novel insights into demographic-specific traffic patterns.

\subsection{Data Acquisition and Model Parameterization}\label{sec:data}
The case study provides a detailed account of the data sources used for the analysis, specifically the traffic volume and speed data. We sourced datasets from the Performance Measurement System (PeMS) Data Source, which is maintained by Caltrans. The datasets consist of two main variables: the total flow $f$, which represents the sum of all vehicles passing by a sensor within a 5-minute period, and the flow-weighted average speed $v$, which represents the average speed of all vehicles passing by a sensor within the same 5-minute period. Less than 5\% of the values were missing; we used linear interpolation to approximate these. Local road information is obtained from the ADMS database, which is managed by the University of Southern California.

We then computed the traffic flow $y_a$ for each link $a$ based on the raw sensor traffic information. We define $S_a = \{s| s\in a\}$ as the set of sensors on link $a$. $|S_a|$ is the number of sensors in the set $S_a$. To determine the length of a link, length($S_a$), we calculated the sum of L2 distances between adjacent sensors along the link. We then used length($S_a$) to calculate the traffic flow. We calculated the traffic flow density of $S_a$ by taking the average traffic flow density of the sensors on link $a$ and multiplying it by the length of the link:
\begin{equation}
y_a = \sum_{s\in S_a}\frac{f_s}{|S_a|v_s}\times \text{length}(S_a)
\label{density}
\end{equation}
In the above equation, $f_s$ represents the traffic flow volume of sensor $s$, while $v_s$ represents the flow-weighted average speed of that sensor. The first term is the flow value of the sensor divided by the average speed of the sensors on the link, multiplied by the total number of sensors on the link. The product of this value and the length of the link gives us the total traffic flow $y_a$ for that link.

In the base model, it is also necessary to define both the DAR and route choice components, as described in section \ref{sec:basic_model}. The DAR can be calculated following the approach detailed in \cite{Ma2018EstimatingData}. For the route choice component, we make the assumption that drivers are equally likely to select between the shortest paths.

The DODE analysis requires constructing the traffic network in LA County and determining the weights of the regularization terms in the optimization model. We examined different network configurations and parameter weights in our analysis to select reasonable values. The optimization problem was solved using stochastic projected gradient descent coded in Python, as introduced in \cite{Ma2018EstimatingData}. All experiments were conducted using High-Performance Computing with Xeon-2640v3 CPU and 64GB of memory.

\FloatBarrier
\subsection{Determining Network Configuration and Objective Value Weights}\label{sec:ne_network}
While the road system in Los Angeles informs the DODE graph, there remain many choices around the network configuration used in our analysis. Additionally, the choice of weights for each regularization term in the objective can impact the final solution. We therefore performed several experiments to determine reasonable choices for these modeling inputs. 

\subsubsection{Network Size}

\begin{figure}[!htb]
        \centering
        \subfloat[Health districts in LA County.\label{sfig:hd}]{\includegraphics[width=0.45\linewidth]{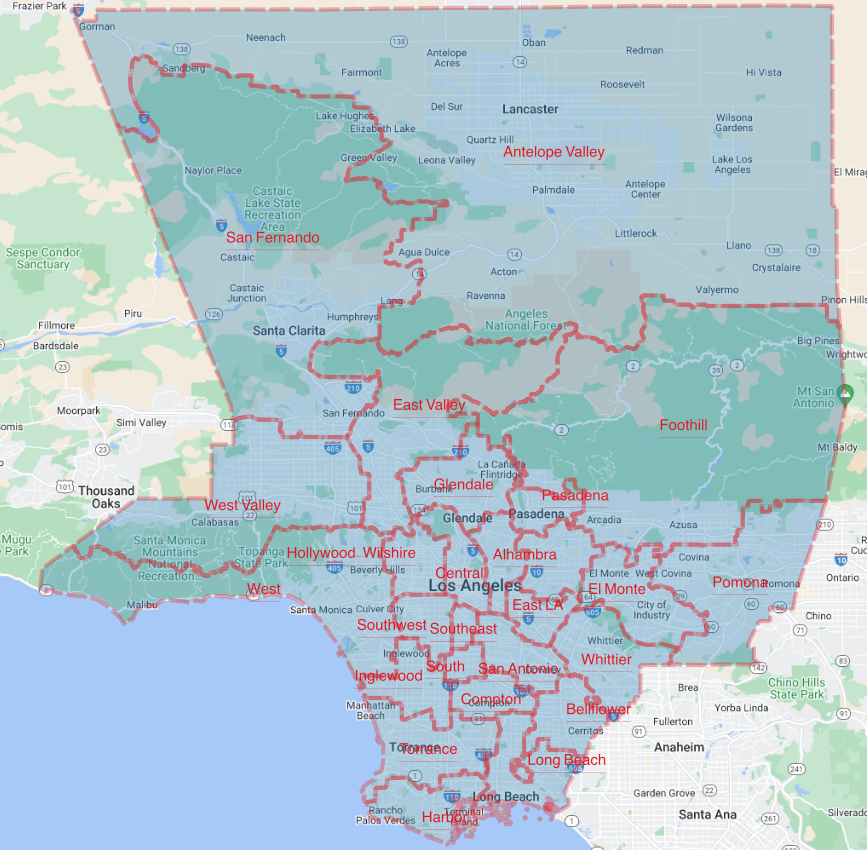}
	}\hfill
	\subfloat[Sensors on highway in LA County\label{sfig:sensor}]
		{\includegraphics[width=0.45\linewidth]{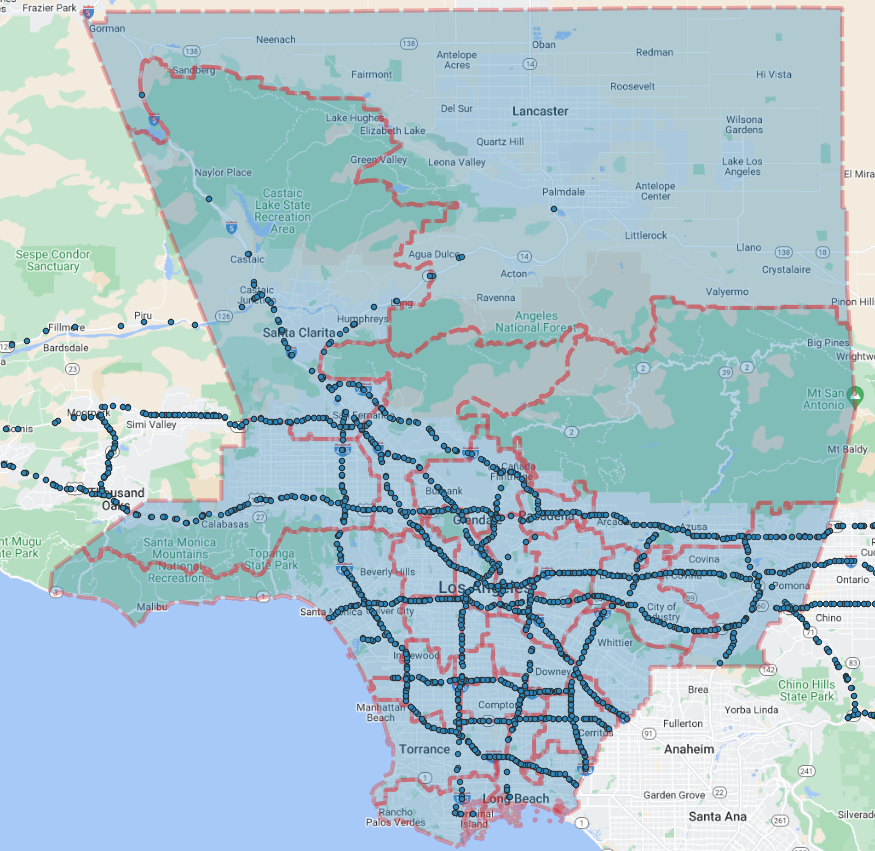}
        }\hfill
	\caption{Health districts and traffic sensors in LA County. Health district boundaries are shown in red. Blue dots indicate sensor locations on highways.}\label{fig:LA_tf}
\end{figure}

Our study relies on the traffic flow data recorded by the sensors on the highway. Figure \ref{sfig:hd} shows the health districts in LA County. There are a total of 26 health districts. We focus on analyzing how people travel between different health districts, and thus we study the aggregate origin and destination flow between each pair of health districts. Figure \ref{sfig:sensor} displays the locations of the highway sensors in LA County. We construct the network along the highways to utilize the sensor information. To construct the traffic network, we first identified potential nodes and edges. In this case, we choose intersections between highways as potential nodes, which allowed us to capture the major points of interchange in the road system. Next, we connected these nodes with links, which corresponded to the physical highways. To increase the resolution of the network, we added additional nodes on long links in urban areas. This ensured that our network accurately represented the actual road system and captured the key points where traffic might change flow or experience congestion. This approach allowed us to construct a detailed model of the traffic network.

\FloatBarrier
\begin{figure}[!htb]
        \centering
        \subfloat[Low\label{sfig:large_network}]{\includegraphics[width=0.32\linewidth]{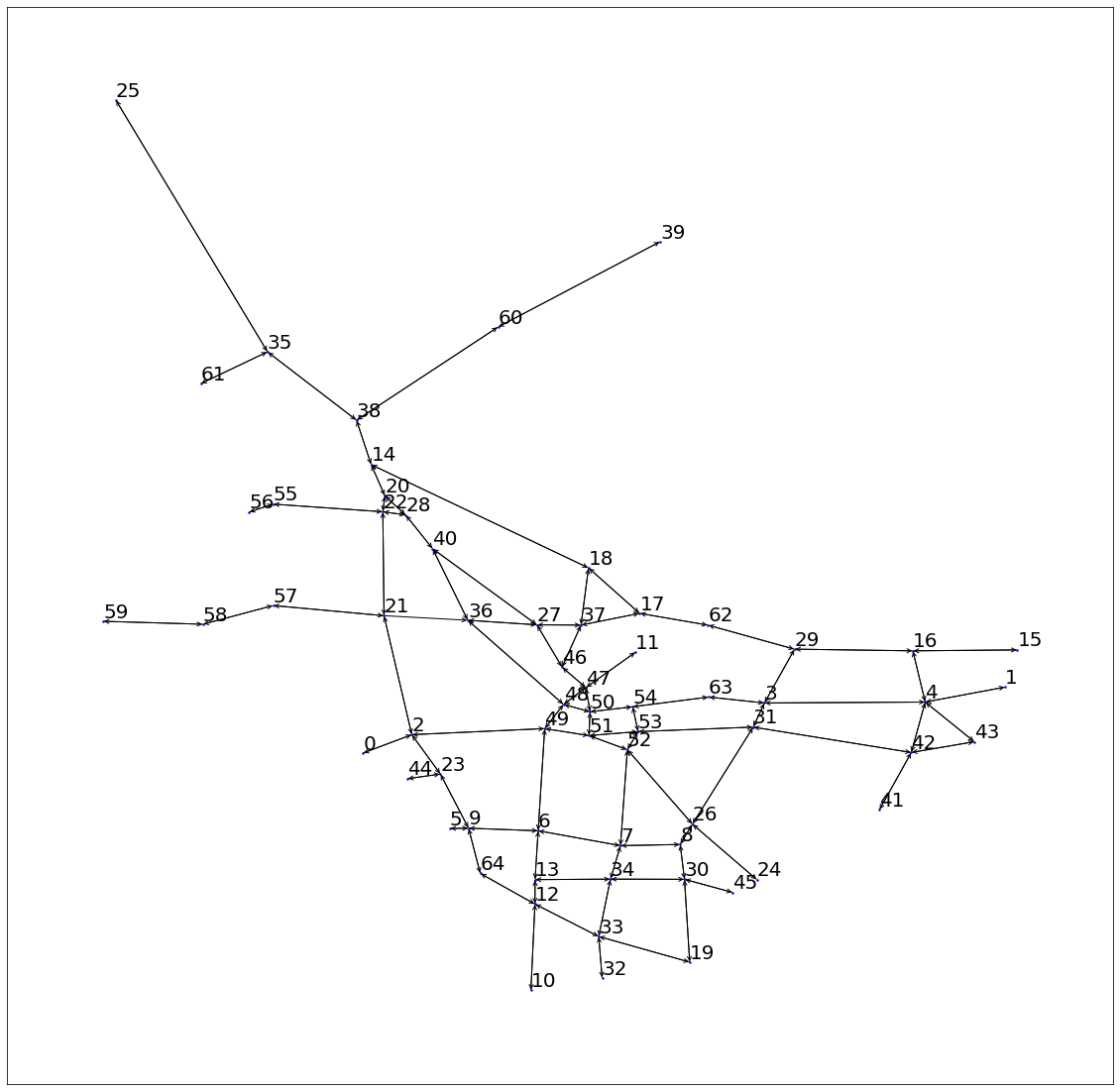}
	}\hfill
	\subfloat[Medium\label{sfig:medium_network}]{
		\includegraphics[width=0.32\linewidth]{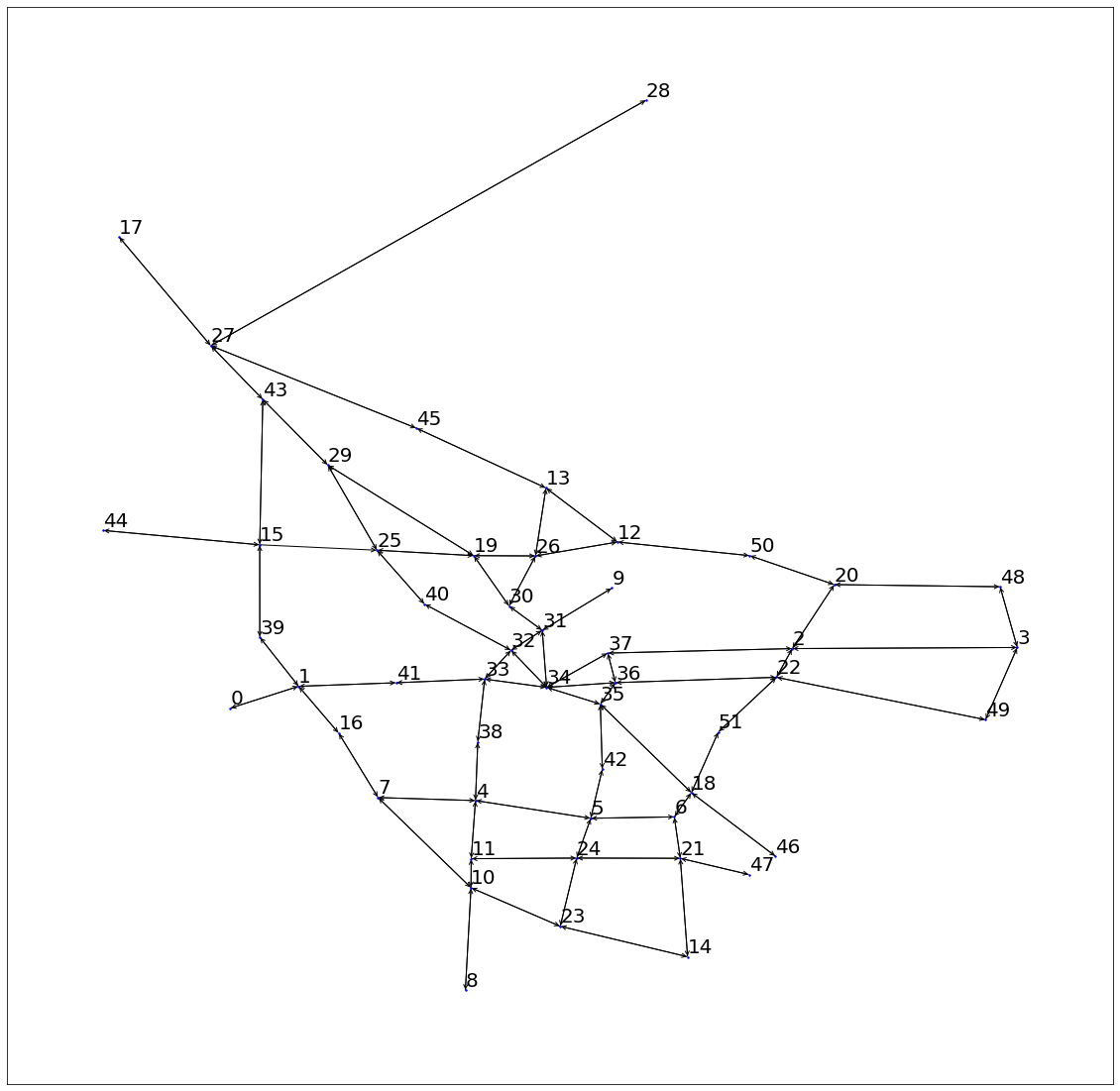}
    }\hfill
        \subfloat[High\label{sfig:small_network}]{
		\includegraphics[width=0.32\linewidth]{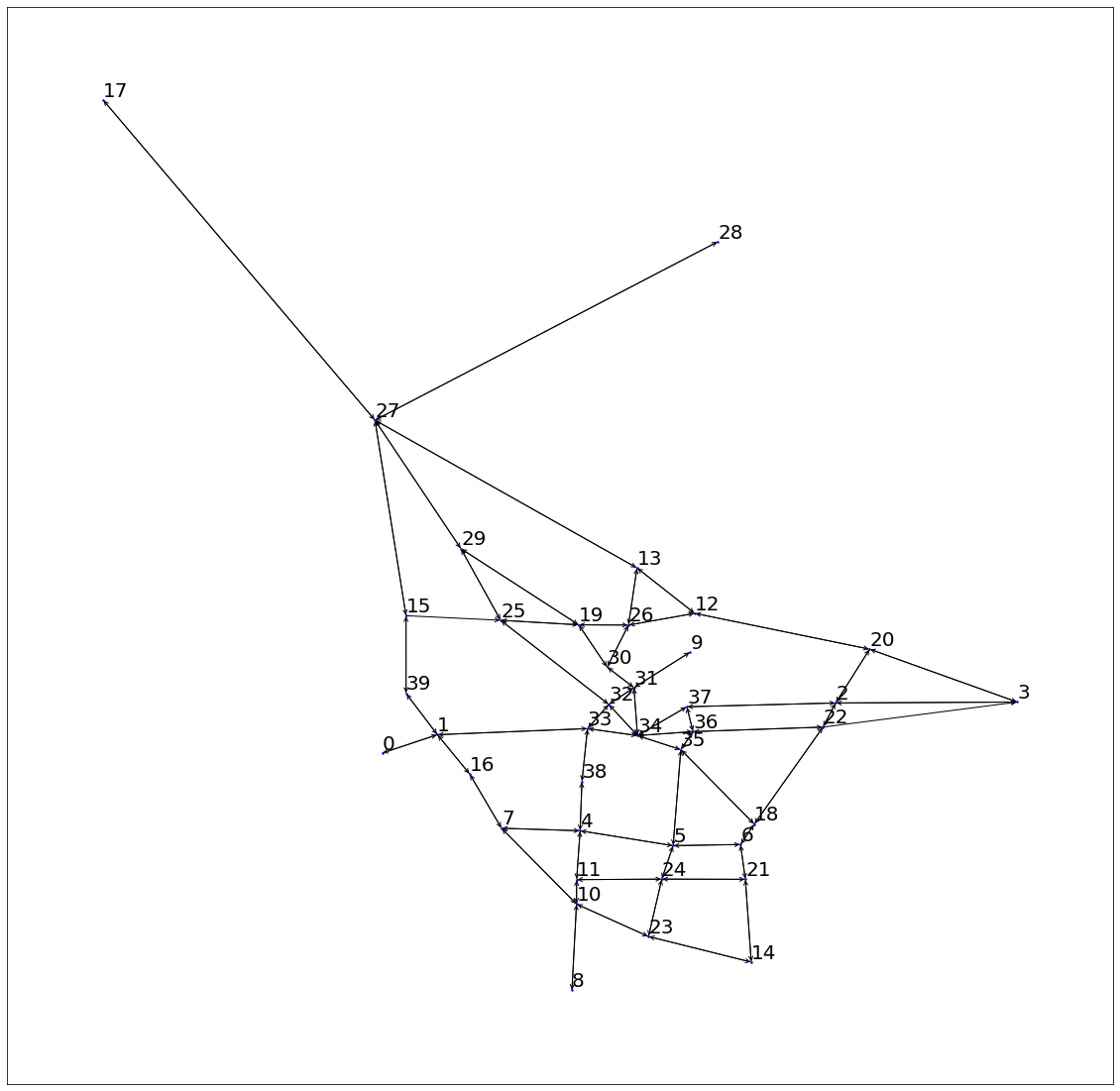}
    }\hfill
    \subfloat[Low\label{sfig:large_network_map}]{\includegraphics[width=0.32\linewidth]{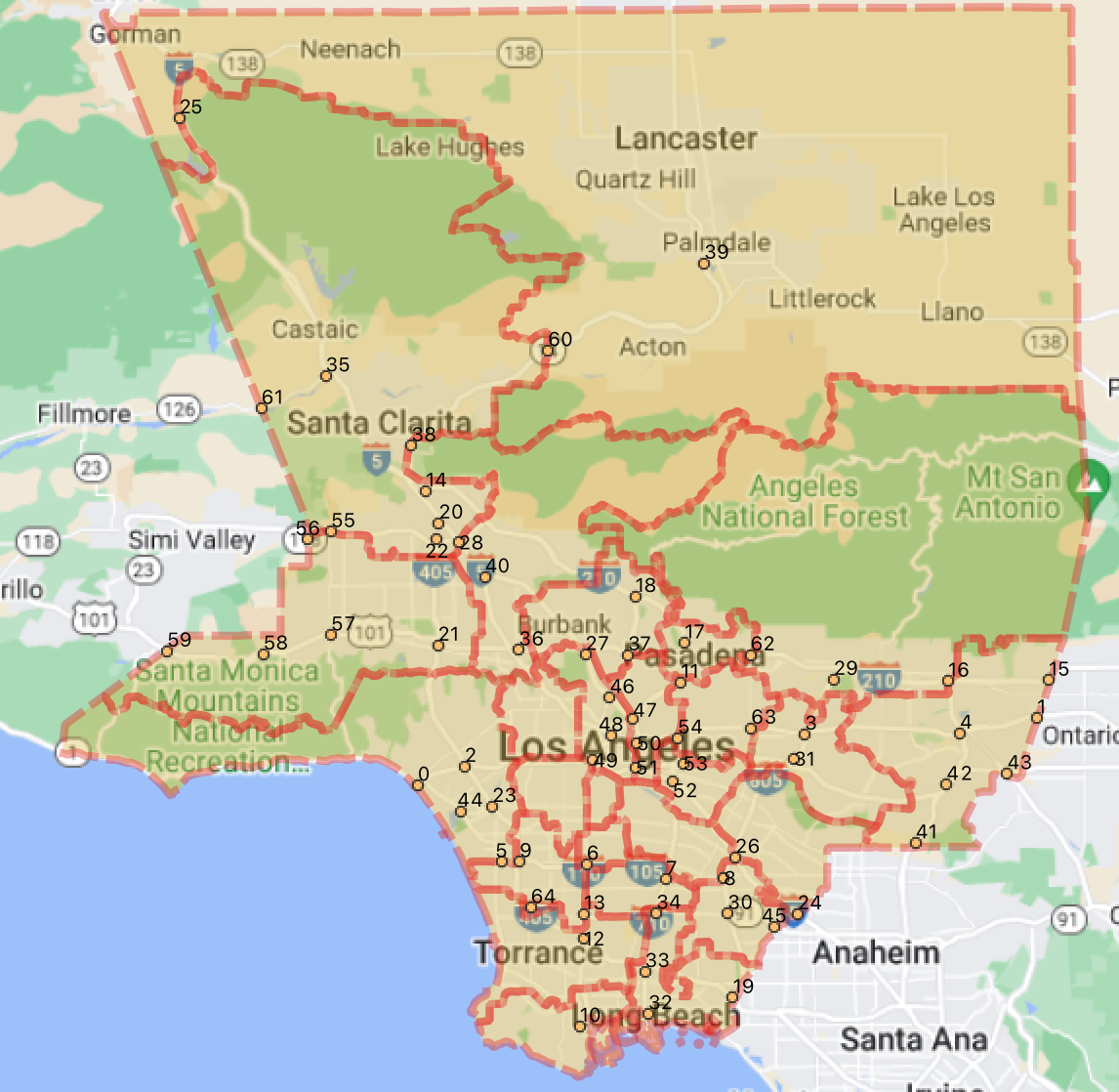}
	}\hfill
	\subfloat[Medium\label{sfig:medium_network_map}]{
		\includegraphics[width=0.32\linewidth]{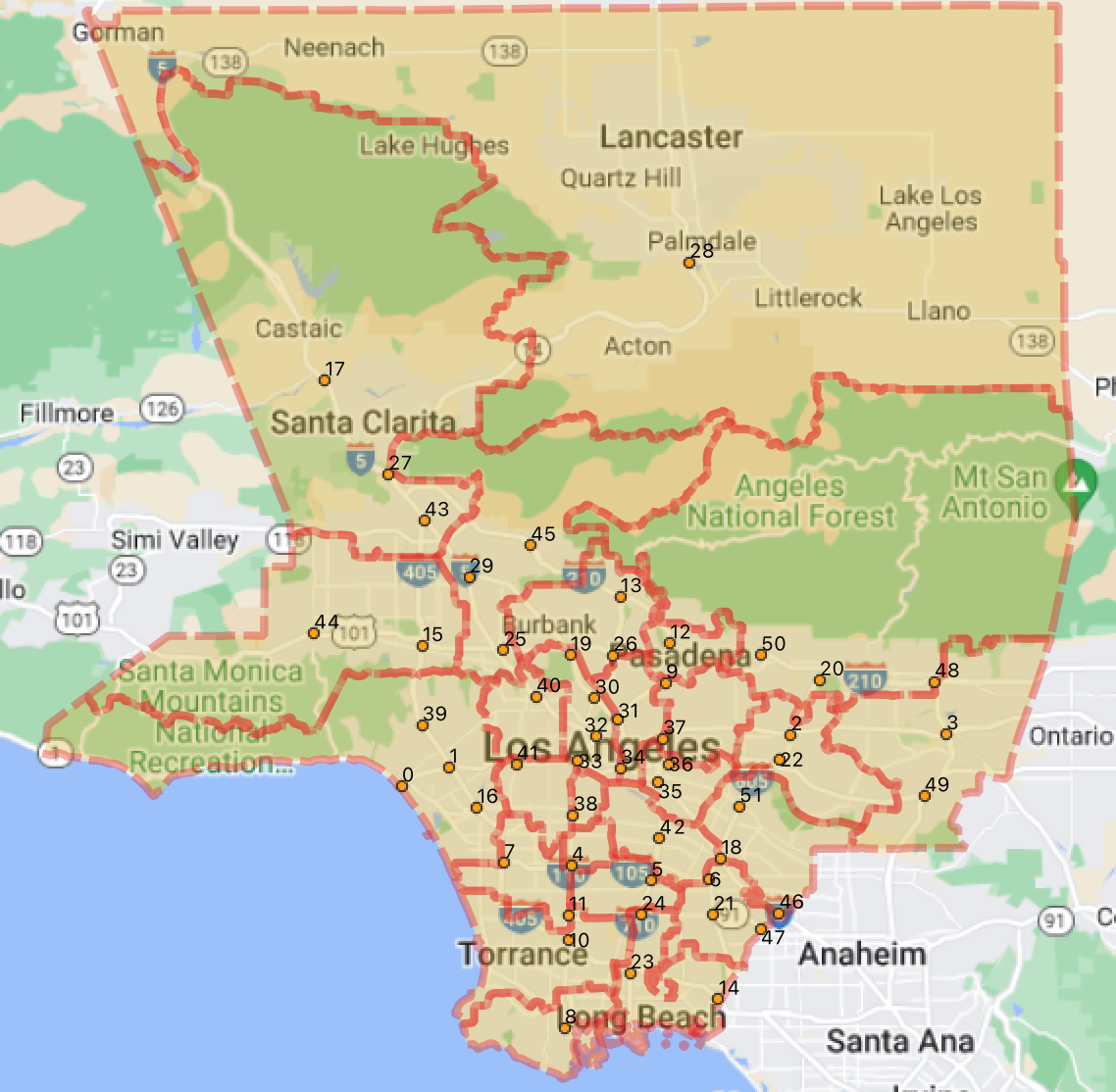}
    }\hfill
        \subfloat[High\label{sfig:small_network_map}]{
		\includegraphics[width=0.32\linewidth]{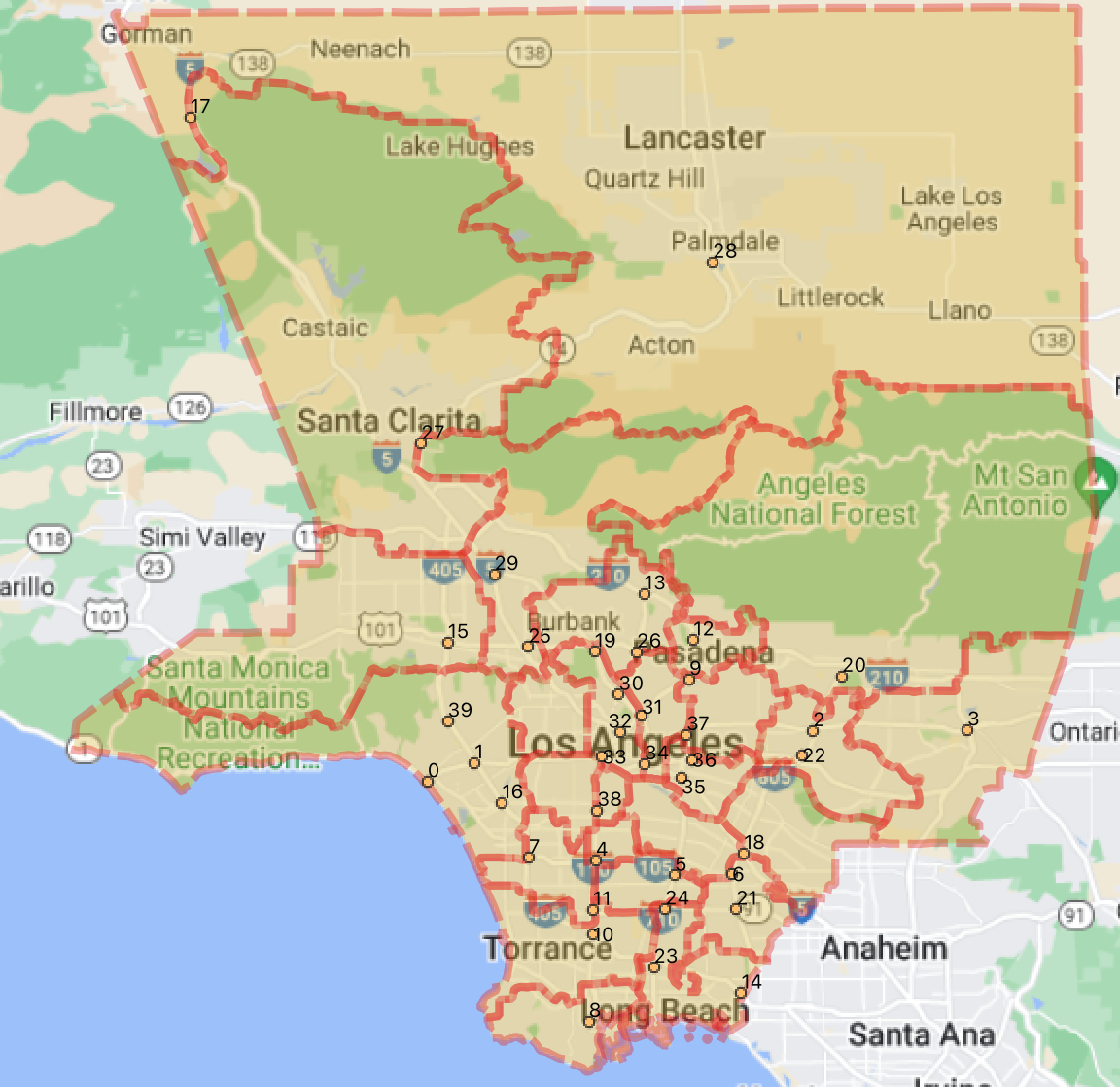}
    }\hfill
	\caption{Different levels of aggregation. Each dot represents both an origin and a destination. The dots are connected by links that follow the highway network in LA County.}\label{fig:LA_nw}
\end{figure}

\begin{table}[hbt!]
\begin{center}
\begin{tabular}{ |c|c|c|c|c| } 
 \hline
 Aggregation level& Number of nodes & Number of ODs & Average computing time (h)& Estimation error ($\epsilon_{b}$)\\ \hline
Low& 65  & 4160 & 12.15 & 1379\\ \hline 
Medium& 52  & 2652 & 2.56 & 1403 \\ 
 \hline
High& 40 & 1560 & 1.37 & 1722\\ 
 \hline
\end{tabular}
\end{center}
\caption{Estimation error and computing time as a function of network size. The average computing time is the average time needed for computing OD flows over a five day period.}
    \label{tab:networks}
     \vspace{-10pt}

\end{table}


However, even with these basic guidelines, there are multiple networks that could be used. We therefore constructed three different networks, each with varying numbers of nodes and edges (see Figure \ref{fig:LA_nw}). The size of the dataset and the computational speed of the optimization model is affected by the number of nodes and edges, and thus we sought to strike a balance between network size and computational efficiency. The first network, depicted in Figure \ref{sfig:large_network}, contained all highway interchanges as nodes, which were connected by the highway network. Additional nodes were added to links that were too long. The network labeled as "low" refers to the network with the highest number of nodes, denoting the lowest level of aggregation. This rule also applies to medium and high aggregation level networks. The medium and high aggregation level networks have fewer nodes and edges but still retain the highway network's structure by applying several rules. Firstly, we reduced the number of nodes in low population density areas. Secondly, in areas where the nodes were densely packed, we eliminated nodes. Finally, in certain rare cases, we may add additional nodes to high population areas if the existing nodes are too sparse. For example, we have added additional nodes 38, 41, and 42 in Figure \ref{sfig:medium_network} as areas around the center of LA County where the density of nodes is insufficient. Then we ran the optimization model on each of these three networks and analyzed their properties, which are displayed in Table \ref{tab:networks}. The table shows that the number of origin and destination (OD) pairs increases exponentially with the number of nodes, which in turn affects the computational time. As the number of OD pairs increases, the computational time grows exponentially.

Table \ref{tab:networks} presents the estimation error $\epsilon_b$ and average computing time for different network aggregation levels. The results indicate that the medium and low aggregation level networks have similar estimation error, but the medium network has a much smaller computing time. Although the high aggregation network runs even faster, it has a much higher estimation error. Thus, we used the medium aggregation level network for our analysis to balance computational time and network complexity. As depicted in Figure \ref{fig:LA_nw}, this network includes at least one node in each health district, with more nodes concentrated in areas of high population density, such as the central, northeast, and west regions. 
\subsubsection{Choice of Weights in Objective Function}
The accuracy of our DODE traffic flow prediction relies on a balanced consideration of link flow evaluation, lower bound determined by local traffic volume, symmetricity, and similarity in total flow outcomes to the classic DODE model. These are captured in the objective function by $\epsilon_{b}$, $\epsilon_{LB}$, $\epsilon_{s}$, and $\epsilon_{\tau}$, respectively. To control their relative importance, the last three have weights ($\eta$, $\beta$, and $\gamma$, respectively). It is important to carefully select these parameter values. To optimize our model and select the most suitable weight values, we assess model outcomes through numerical evaluation over a range of weights. Since $\epsilon_{LB}$ serves as a lower bound for total nearby arterial traffic flow, it is always fulfilled due to the inclusion of slack variables. Since we only require that the slack variables be non-negative, we do not concern ourselves with their magnitude. We can thus set the weight parameter $\eta$ arbitrarily; in our numerical examples, we use a value of 1.

In the next step, we conduct experiments using various combinations of $\beta$, $\eta$, and $\gamma$. Our goal was to strike a balance between the three sources of error: estimation error, symmetric error, and total flow. Our approach involves first determining the total traffic flow by running the base optimization model with $\eta = 0$, $\beta = 0$, and $\gamma = 0$. Once we established this baseline value, we introduced local road traffic flow lower bound, symmetricity, and total flow errors into the model and tested it with different weights of error terms. To ensure a comprehensive analysis, we did this by sweeping over a range of $\beta $ and $\gamma$ values. The results are presented in Table \ref{tab:error_parameter}. 

\FloatBarrier
\begin{table}[H]
\begin{center}
\centering
\begin{tabular}{ |c|>{\centering\arraybackslash} p{1in}|>{\centering\arraybackslash} p{1in}|>{\centering\arraybackslash} p{1in}|>{\centering\arraybackslash} p{0.8in}| } 
 \hline
 Parameter values &Estimation Error \newline ($\epsilon_{b}$) & Symmetric Error \newline ($\epsilon_{s}$) & Local road traffic \newline ($\epsilon_{LB}$)  & Total Flow \newline (vehicles)\\ \hline
$\beta = 0$, $\eta = 0$, $\gamma = 0$& 1403 & 599 & 4835 & 1966014\\ \hline
$\beta = 1$, $\eta = 1$, $\gamma = 1$ & 1466 & 414&  2318 & 1963851 \\ 
 \hline
$\beta = 10$, $\eta = 1$, $\gamma = 1$& 1490 & 104& 2297 & 1965013\\ 
 \hline
$\beta = 100$, $\eta = 1$, $\gamma = 1$& 1535 & 110& 3008 &1959267 \\ \hline
$\beta = 1e3$, $\eta = 1$, $\gamma = 10$ & 1710 & 111 & 4420 &1953910\\ 
 \hline
$\beta = 1e5$, $\eta = 1$, $\gamma = 10$& 1724 & 113& 4544 &1952400\\ 
 \hline
$\beta = 1e20$, $\eta = 1$, $\gamma = 100$& 2406 & 18& 1598 &1163434\\ 
 \hline
\end{tabular}
\end{center}
\caption{Errors under different parameter value selections.}
    \label{tab:error_parameter}
\end{table}

Table \ref{tab:error_parameter} displays the estimation error, symmetric error, local road traffic error, and total flow values obtained from the optimization model under different parameter settings. We see that non-zero values of $\beta$ and $\eta$ significantly improve the symmetric error. When $\beta$, $\eta$, and $\gamma$ are increased from 0 to 1, the symmetric error is reduced by 30\% compared to the original symmetric error. Additionally, the local road traffic error decreases by more than 50\%, while the estimation error only increases by around 6\%. When the value of $\beta$ is increased to 10, we observe an 82\% reduction in the symmetric error. Furthermore, the error related to local road traffic and the estimation error remains relatively stable. When we further increase the weight of the symmetric error by setting $\beta$ in the range of 10 to 1e5, the decrease in the symmetric error is not substantial. However, other errors do increase. Moreover, the total traffic flow decreases as $\beta$ increases, but it remains within a reasonable range. If we raise $\beta$ to an extremely large value, such as $1e20$, the symmetric error approaches zero, and there is also a decrease in local road traffic. However, this improvement comes at the cost of a substantial increase in the estimation error.

To restrict the decrease in total traffic flow, we tune the $\gamma$ parameter accordingly. As a result, the total traffic flow does not show much difference from the original total traffic flow. However, if we set the weight of the symmetric error to a very large number (e.g., 1e20), the symmetric error approaches zero, but the total traffic flow decreases by around 41\%, and the estimation error increases by 71\%.

We need to find a balance between the estimation error, symmetric error, and total traffic flow. Therefore, we choose a pair of $\beta$ and $\gamma$ values that can achieve this balance. From Table \ref{tab:error_parameter}, we can see that the pair $\beta = 10$ and $\gamma = 1$ reduces the symmetric error by 82\%, while the estimation error and total flow remain at a similar level as the base model. Therefore, we select $\beta = 10$ and $\gamma = 1$ for the remaining analysis.
\FloatBarrier
\subsection{DODE Estimation Outcomes Capture Epidemic Trends in Traffic Flow Patterns}
With the network configuration and parameters chosen, we can use the optimization model to estimate DODE outcomes. We do this for the LA County data from 2019 and 2020. In 2020, LA County imposed stay-at-home order due to the COVID-19 pandemic, severely changing traffic flow patterns. We can therefore understand these patterns through comparison of 2019 and 2020 DODE outcomes from our model. 
\begin{figure}[H]
        \centering
        \subfloat[Box plot of traffic flow in 2019 and 2020, with outliers represented by dots. The box displays the interquartile range of OD traffic flow across all districts, with the median indicated by a short horizontal line. The highest and lowest short horizontal lines outside the box indicate the maximum and minimum values, respectively. The x-axis of the graph represents the timeline and includes significant time spots when COVID-19 policies and interventions were implemented. \label{sfig:tf_boxplot}]{\includegraphics[width=1\linewidth]{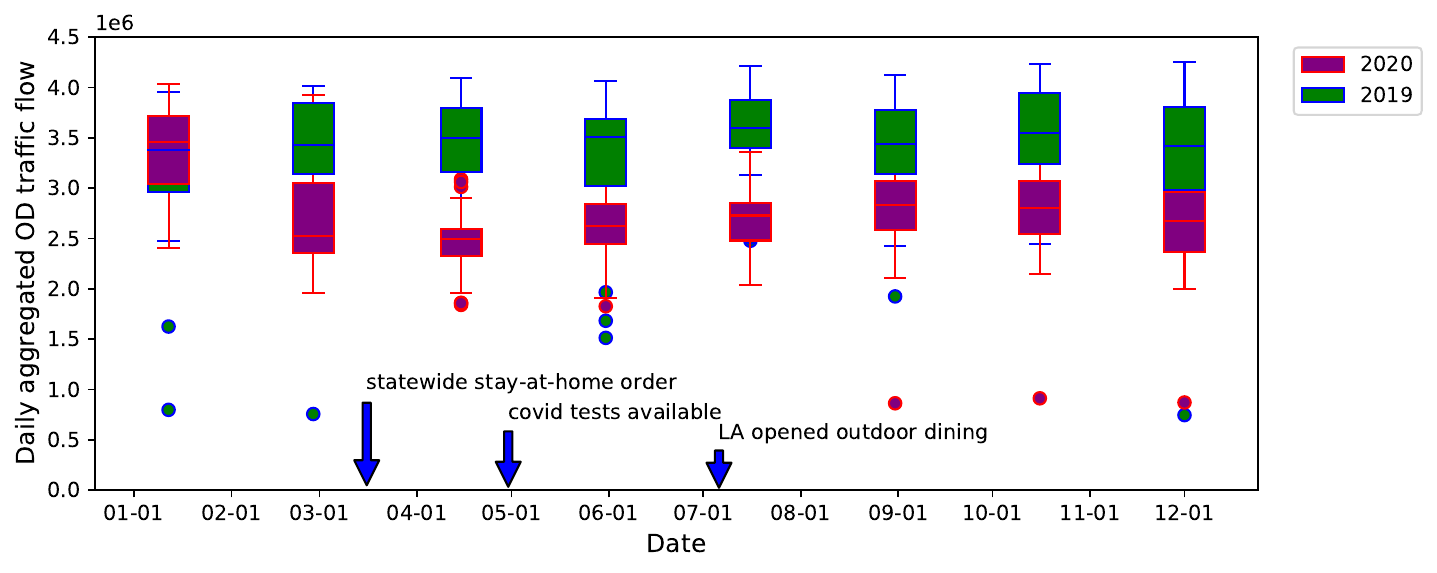}
	}\hfill
	\subfloat[Aggregated daily traffic flow on 03/20/2019 and 03/20/2020. \label{sfig:tf_daily}]{\includegraphics[width=0.5\linewidth]{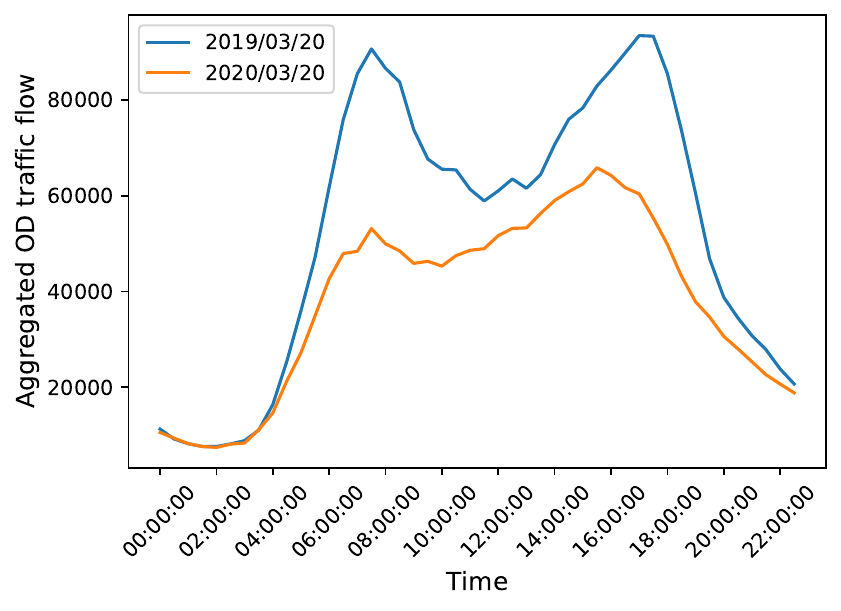}
    }\hfill
	\caption{Traffic flow comparison in 2019 and 2020}\label{fig:tf_compare}
\end{figure}

In particular, we analyze two outcomes: the change in the DODE model traffic volume over 46-day intervals and the OD daily traffic flow within 24 hours. To compare the traffic volume in 2019 and 2020, we divide the two-year period from January 12, 2019 to December 31, 2020 into 16 intervals. Every year is divided into 8 intervals, with each interval comprising 46 days. By conducting this comparative analysis, we aim to gain a comprehensive understanding of the effects of the COVID-19 pandemic on traffic flow and commuting behavior in LA County.

We calculate the total demand for all OD pairs by summing the traffic volume of each pair for every day and compare intervals with similar starting and ending dates in 2019 and 2020 in Figure \ref{sfig:tf_boxplot}. This figure shows that aggregate OD demands where generally lower in 2020 compared to 2019 after February. 
This could be because individuals were concerned about the rapidly rising COVID-19 case counts at that time, and traffic flow volumes remain depressed through the rest of the year. The traffic flow distribution across all intervals in 2019 remained similar. During the first period, the OD demand distributions for 2019 and 2020 had significant overlap with each other, indicating that there was no statistically significant difference in traffic flow between 2019 and 2020 before the stay-at-home order (see Figure \ref{sfig:tf_boxplot}).  However, starting from the second period, traffic flow in 2020 was depressed relative to 2019 levels, reaching its lowest point in the third period, which was just after the start of the stay-at-home order in Los Angeles. Traffic flow remained at a low level for the following two periods and increased slightly in the last three periods in 2020. This suggests that traffic flow patterns in 2020 were indeed greatly reduced due to the stay-at-home order. This demonstrates that our DODE model can capture general traffic trends as suggested by other sources that showed decreased traffic volume during this period \citep{Abdullah2020ExploringPreferences}.

However, our data can provide additional information besides these general trends. In addition to analyzing DODE model traffic demand across different time periods, we also examined the daily OD traffic flow. Figure \ref{sfig:tf_daily} depicts the total traffic flow from the DODE model on an example day (March 20 in 2019 and 2020). There are two traffic peaks over the day, occurring around 7am and 5pm, which aligns with the typical work schedule. These peaks indicate the rush hours when a significant number of individuals are traveling to and from their workplaces. However, the overall traffic flow was lower in 2020, as expected due to the stay-at-home order. While there was not much difference in traffic flow during non-peak hours, the gap between 2019 and 2020 became much larger during peak hours. This reduction in traffic flow clearly indicates that the pandemic had a substantial impact on commuting travel activity during the stay-at-home order period.

These preliminary analyses demonstrate that the model can capture reasonable expected trends in 2019 and 2020 data. It is evident that there has been a significant reduction in traffic flow from origin to destination in both periodical and daily traffic demands during the COVID-19 pandemic in comparison to 2019. However, the true power of the DODE model lies in analysis of origin and destination estimations, which we examine in the subsequent section. 
\FloatBarrier
\subsection{OD Flow between Health Districts}
This DODE model can be used as a key tool for understanding human behavior in pandemic environments, which can provide public health officials with critical information about how health policies like stay-at-home orders impact population flows in the county. By analyzing the patterns of traffic movement across health districts, researchers can identify how public health policies deferentially affect different subgroups in the population, whether intentionally or not. A detailed understanding of traffic flow can also help in the allocation of healthcare resources, such as the distribution of medical supplies or the placement of healthcare facilities. 

We use our DODE model to identify traffic flow between health districts. We first analyze the traffic flow between health districts during each time period separately, treating the estimation of traffic flow on each date as a data point. To determine whether there is a significant difference between the OD flows in 2019 and 2020, we aggregate the daily traffic flow for each OD pair across each interval. Then, we conduct a paired t-test on traffic flow values in 2020 and 2019 within each interval for all OD pairs between health districts with sufficient data (we do not consider OD pairs with traffic flow values less than 1000 vehicles per day).

For each time period, Table \ref{tab:t_test} reports the total number of OD pairs for which traffic flow increased or decreased in 2020 compared to 2019. The table also shows the number of OD pairs for which the increase or decrease in traffic flow was statistically significant at the 0.05 level during the specified time period, indicated by the p-value. The table reveals that the number of OD pairs with increased traffic flow got lower as time progressed, whereas the number of OD pairs with decreased traffic flow increased. This trend became more pronounced after the stay-at-home order was issued. In the first period, prior to the stay-at-home order, there are a similar number of OD pairs that had increased and decreased traffic flow (101 and 108), suggesting that there was not a time trend between 2019 and 2020 traffic flow, and only around $25\%$ of OD pairs are significant. However, during the period from February 27th (just before the stay-at-home order started on March 15th) to April 14th, the total number of increased OD pairs decreased from 101 to 54, and we see fewer OD pairs with statistically significantly increased traffic flow (12 OD pairs) and more with decreased traffic flow (65 OD pairs). This trend persists for the remainder of the year. After the period of April 14th, the majority of OD pairs continued to experience a decline in traffic flow. On the period August 30, the number of statistically significantly decreased OD pairs reduced from 104 to 67. This suggests that traffic flow started to return to close to normal levels as the impact of the pandemic on travel reduced. However, it is worth noting that even by the end of the year, the number of decreased OD pairs (160 OD pairs) is still much higher than the number of statistically significant increased OD pairs (43 OD pairs), indicating that the pandemic still had an effect on traffic flow at that time. These findings demonstrate that the expected trend of reduced traffic flow due to the stay-at-home order can be clearly seen even in our OD outcomes.
\FloatBarrier
\begin{table}[H]
\begin{center}
\begin{tabular}{ |p{2cm}|p{2cm}|p{3cm}|p{3cm}|p{3cm}| } 
 \hline
Time period & Total number of increased OD pairs & Number of significant increased OD pairs (p-value $\leq 0.05$) & Total number of decreased OD pairs & Number of significant decreased OD pairs (p-value $\leq 0.05$)\\\hline
01-12 & 101 & 30 & 108 & 20\\\hline
02-27 (Stay-at-home order start at 03-15) & 54 & 12 & 150 & 65\\\hline
04-14 & 40 & 7 & 165 & 104\\\hline
05-30 & 50 & 8 & 150 & 71\\\hline
07-15 & 39 & 12 & 154 & 89\\\hline
08-30 & 66 & 11 & 147 & 67\\\hline
10-15 & 40 & 6 & 166 & 89\\\hline
11-30 & 43 & 6 & 160 & 87\\\hline
\end{tabular}
\end{center}
\caption{t-test on traffic flow between 2019 and 2020. It summarizes the results of a t-test analysis on traffic flow between 2019 and 2020 during different time periods. The table includes the following columns: the period of dates for which the t-test was performed, specified by the start date of the period, the total number of OD pairs for which traffic flow increased/decreased in 2020 compared to 2019 during the specified time period, the number of OD pairs for which the increase/decrease in traffic flow was statistically significant at the 0.05 level during the specified time period.}
    \label{tab:t_test}
\end{table}

\begin{figure}[htbp]
        \centering
        \includegraphics[width=1\linewidth]{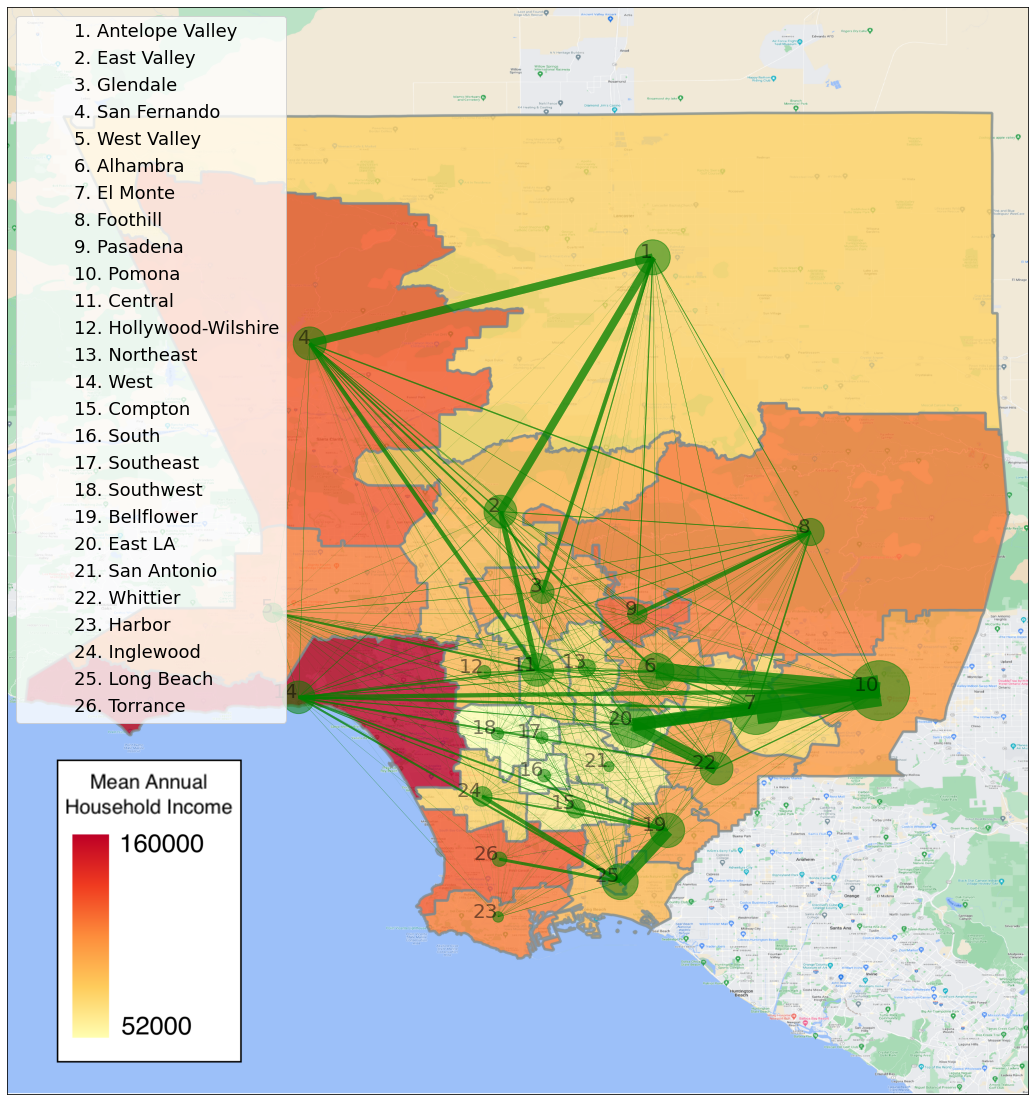}
	\caption{Mean traffic flow in 2020-04-14 to 2020-05-29 period. The green edge indicates the amount of mean traffic flow. Larger edge widths and diameters indicate higher traffic flow. The size of the circle represents the amount of aggregated traffic inflow and outflow of each district. The health districts are represented by different colors based on the average annual income of families residing within the district.}\label{fig:hd}
\end{figure}

However, beyond these general traffic flow findings, our DODE approach allows us to understand traffic flow by OD pair. Figure \ref{fig:hd} depicts the relative traffic flow between health districts in 2020 where the width of each edge scales with traffic flow and the radius of the circle indicates the relative amount of total inflow and outflow.  Our analysis shows that there was a significant amount of traffic flow between residential districts located on the east of LA county, such as Alhambra (district 6), El Monte (district 7), Pomona (district 10), Bellflower (district 19), Whitter (district 22) and Long Beach (district 25). These districts are characterized by high population densities, which contribute to the higher traffic flow observed in these areas. Pomona (district 10) accounts for the highest total traffic flow. This can be attributed to its large population density and its role as a connecting point to the eastern side of LA County. Given its location, individuals traveling to the east of LA County often pass through Pomona district. We also speculate that this trend may be because of traffic from households with homes in the east who had commuted to downtown LA. Additionally, we observed that there was a considerable amount of traffic traveling from downtown districts, such as Central (district 11), to suburban areas like San Fernando (district 4) and East Valley (district 2), indicates a notable commuter movement between the bustling downtown area and the relatively quieter suburban areas. Furthermore, San Fernando and East Valley districts serve as important connectors to the northern areas of Los Angeles County. This connectivity contributes to the flow of traffic not only between downtown and the suburbs but also towards the northern regions of the county. Our analysis using the DODE approach uncovers significant traffic flow patterns between residential districts, and reveals the dynamics of traffic movement in LA County.

Figure \ref{fig:hd} also displays the income distribution within each health district. We see that there is considerable heterogeneity in traffic flow patterns by OD. Given this observation, we aim to delve deeper into the discrepancies and disparities in the effects of the stay-at-home order on OD traffic. We speculate that during the stay-at-home order, people may have reduced non-essential travel and not needed to commute for work as businesses shuttered or converted to work-from-home policies. In the subsequent Section \ref{sec:income}, we provide a more in-depth analysis and discussion regarding the relationship between traffic flow and income. By exploring this relationship, we can gain a better understanding of how income levels influence traffic flow patterns and vice versa.
 
\FloatBarrier
\subsection{DODE Estimation Outcomes Can Uncover Patterns by Demographic Characteristics}\label{sec:income}

Before the COVID-19 pandemic, people typically traveled to work regularly, but the stay-at-home order encouraged working from home and reductions in leisure travel. However, the availability of working from home may differ by industry and type of work --- essential workers continued to travel to work, and high income office workers were more likely to have viable work from home option, as shown in the work of \cite{Jay2020NeighbourhoodStates}. However, many businesses also shuttered during the stay-at-home order, and this led to layoffs that disproportionately affected lower income workers \citep{Bartik2020TheExpectations}. There are therefore several reasons why income may be correlated to traffic flow pattern changes during the stay-at-home order. 

\begin{figure}[H]
        \centering
\includegraphics[width=0.5\linewidth]{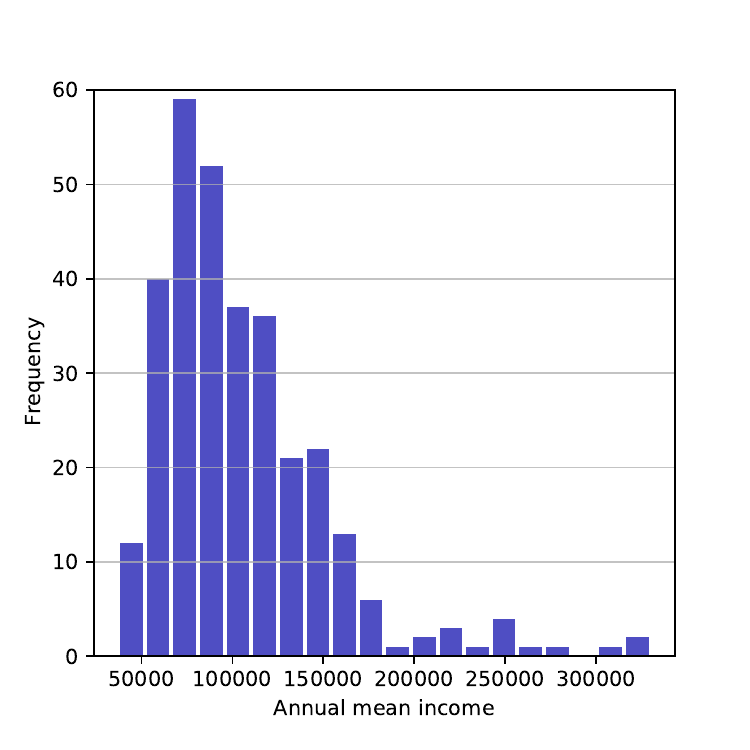}
	\caption{Histogram of Average Annual Household Income by Zipcode. The maximum mean income is more than 330,000 and the minimum mean income is below 50,000.
}\label{fig:income_zip}
\vspace{-10pt}
\end{figure}

To calculate average annual health district incomes for LA County, we used data from the \cite{U.S.CensusBureauLosBureau} by zipcode in LA County. Note that the income distribution is right skewed, with a long right tail (see Figure \ref{fig:income_zip}). Each health district $h$ is comprised of several zipcodes, which we denoted by ${z_i\in h}$. Since health districts can cut across zipcodes, we considered the population $population_i$ of each zipcode $z_i$ within the health district and the proportion of the geographical area of zipcode $z_i$ that falls within health district $h$ when calculating the average district income (equation \eqref{eq:income}). The district income values are shown in Figure \ref{fig:hd}, where we observe that suburban areas of LA tend to have higher average incomes than the centrally located downtown districts, as expected. The average district income values vary from \$52,217 to \$156,975, with a county average of \$90,847. 
\begin{equation}
    income_h = \displaystyle\sum_{z_i\in h} \frac{population_i}{\displaystyle\sum_{z_i\in h}population_i}\frac{area(h\cap z_i)}{area(z_i)} income_i
\label{eq:income}
\end{equation}

Figure \ref{fig:hd} visualizes the annual mean income of each health district overlaid with changes in traffic flow from 2019 to 2020. The plot illustrates that a majority of the OD pairs experienced a decrease in traffic flow. Notably, the outer areas of LA County with higher incomes appear to be more affected by the stay-at-home order compared to the centrally located areas. Moreover, the data suggests that districts with higher mean income, such as the West (district 14), have experienced a relatively greater reduction in traffic flow compared to their original levels. By contrast, lower income districts like South (district 16), Southeast (district 17), and Southwest (district 18) have experienced a relatively smaller decrease in traffic flow. This figure shows that there may be an underlying relationship between changes in traffic and income within LA County during the stay-at-home order, so we therefore investigate the relationship between traffic reduction and income in LA County.

To do this, we wish to compare the percentage difference in traffic flow and a measure of average annual income across OD pairs. For each OD pair $(R_i, R_j)$, the percentage difference in traffic flow is given in Equation \eqref{eq:percentage}, where $f^{year}_{R_i,R_j}(t)$ represents the average traffic flow from health district $R_i$ to health district $R_j$ during time period $t$ in the specified year.
\begin{equation}
    \text{Percentage change in traffic flow} = \frac{f^{2020}_{R_i,R_j}(t) - f^{2019}_{R_i,R_j}(t)}{f^{2019}_{R_i,R_j}(t)}
    \label{eq:percentage}
\end{equation}

We next turn to calculate an income metric for each OD pair. However, we need to consider the income levels at both the origin and destination locations. One natural metric for income to use would be the average district income between the two districts in each OD pair. However, because it is possible that the origin and destination districts may have vastly different incomes, this metric it might obscure underlying patterns. It may also not be useful to only examine income at the origin or at the destination, as we assume traffic flow is symmetric and this would not distinguish between traffic flows to high or low income destinations from similar income origins (or vice versa). We therefore instead use the minimum and maximum in average annual district income. The maximum calculated using $\max\left(income_{R_i}, income_{R_j}\right)$ for each OD pair $(R_i, R_j)$. Similarly, the minimum is given by $\min\left(income_{R_i}, income_{R_j}\right)$.
\begin{figure}[!htbp]
        \centering
        \subfloat[Percentage change in traffic flow between 2019 and 2020 with maximum on income between health districts.\label{sfig:income_max_scatter}]{\includegraphics[width=0.7\linewidth]{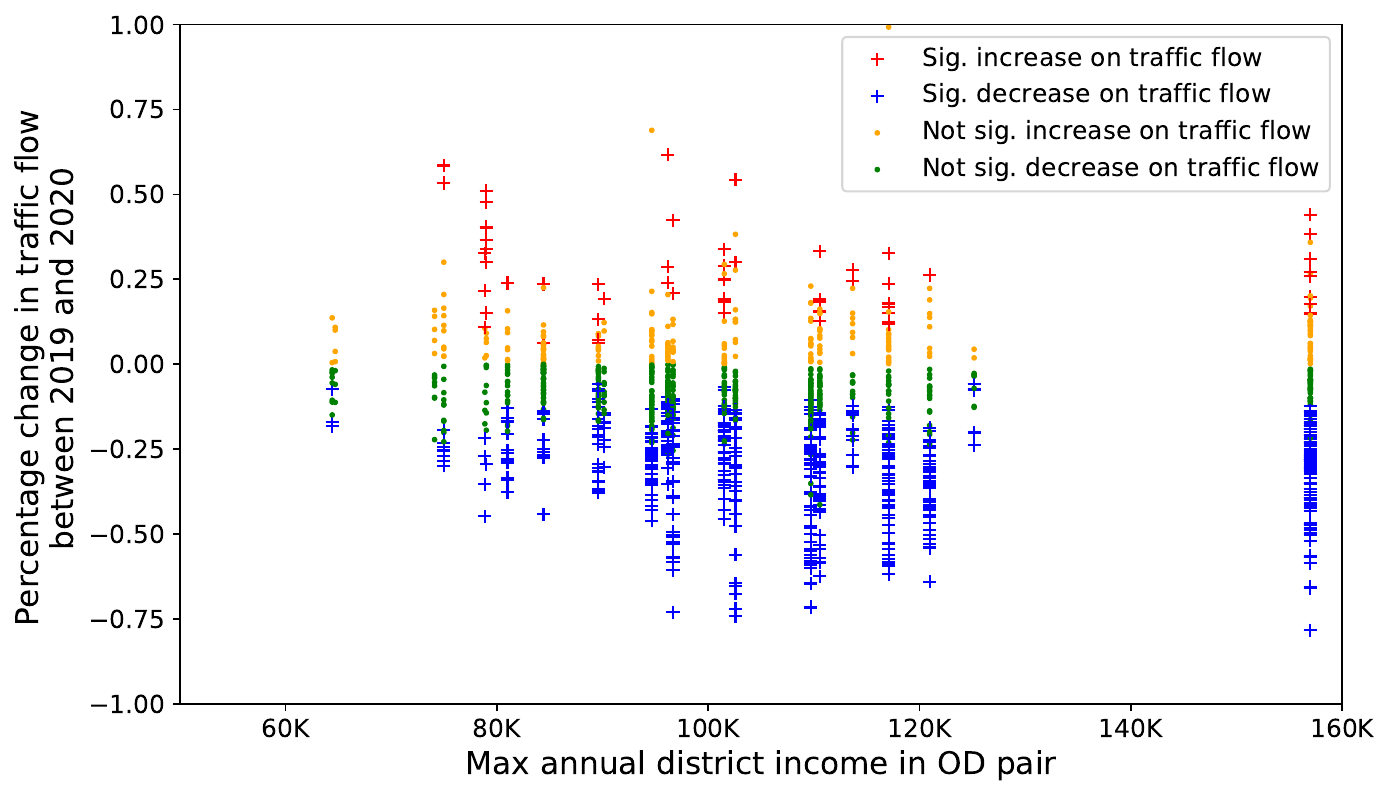}}
	\hfill
        \subfloat[Percentage change in traffic flow between 2019 and 2020 with minimum on income between health districts\label{sfig:income_min_scatter}]{\includegraphics[width=0.7\linewidth]{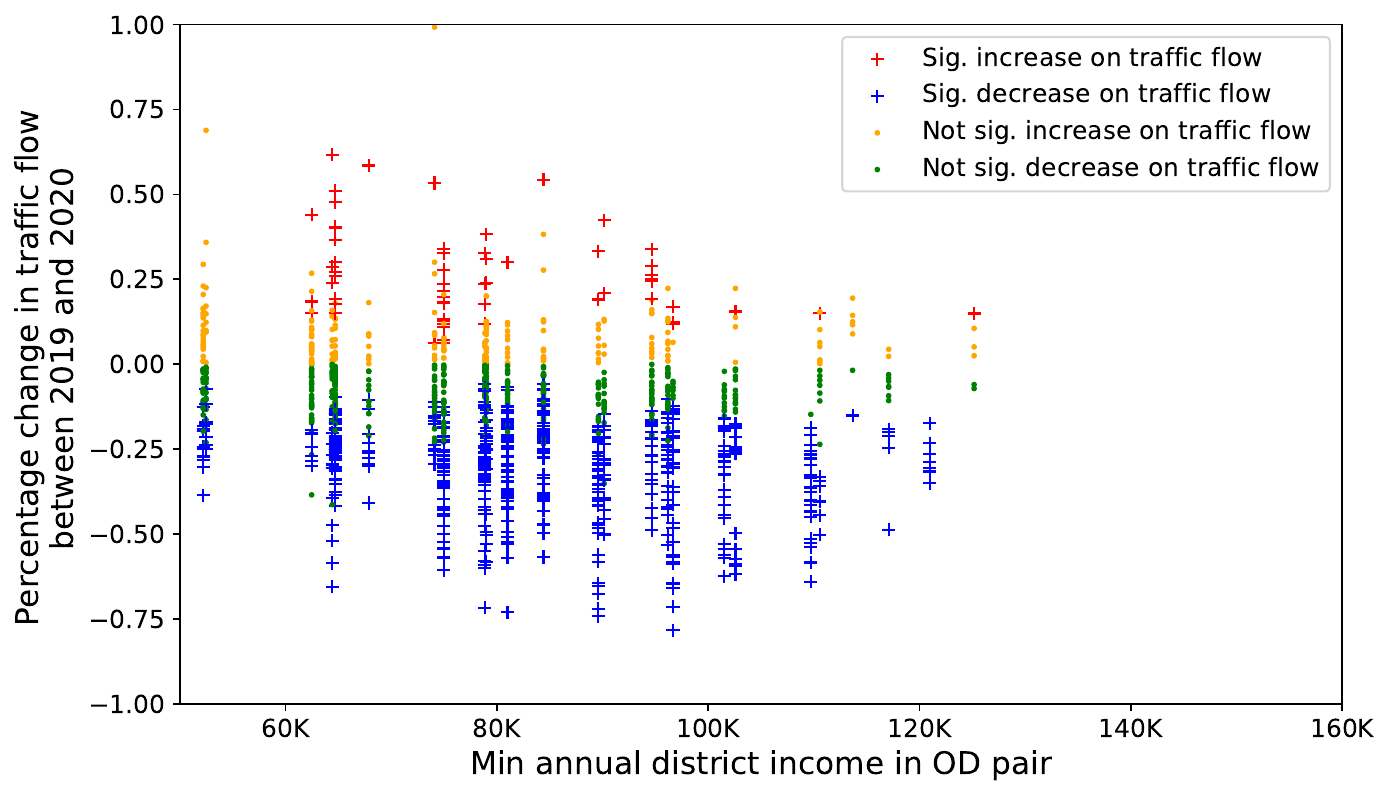}}
	\caption{Percentage change in traffic flow between 2019 and 2020. The figure displays each dot as an average OD pair's traffic flow difference between 2019 and 2020 in a given time period. blue and red plus signs represent OD pairs with statistically significant decreases and increases, respectively. The green and orange dots indicate OD pairs have non-significant decreases and increases.}\label{fig:income_pert}
\end{figure}

Figure \ref{fig:income_pert} visualizes the relationship between the percentage change in traffic flow between 2019 and 2020 by the maximum and minimum incomes in each OD pair. In the figure, blue and red plus signs represent OD pairs with statistically significant decreases and increases, respectively. The green and orange dots indicate OD pairs have non-significant decreases and increases. From the figure, there are more OD pairs with statistically significantly decreased traffic flow than with increased traffic flow, as is expected due to the stay-at-home policy in 2020. 

\begin{figure}[!htbp]
        \centering
        \subfloat[Density of OD pairs by \emph{maximum} income in each OD pair. The dotted line represents the income at \$80,000, which corresponds to the 38th percentile of the income distribution in the health district. \label{sfig:income_max_kde}]{
		\includegraphics[width=0.7\linewidth]{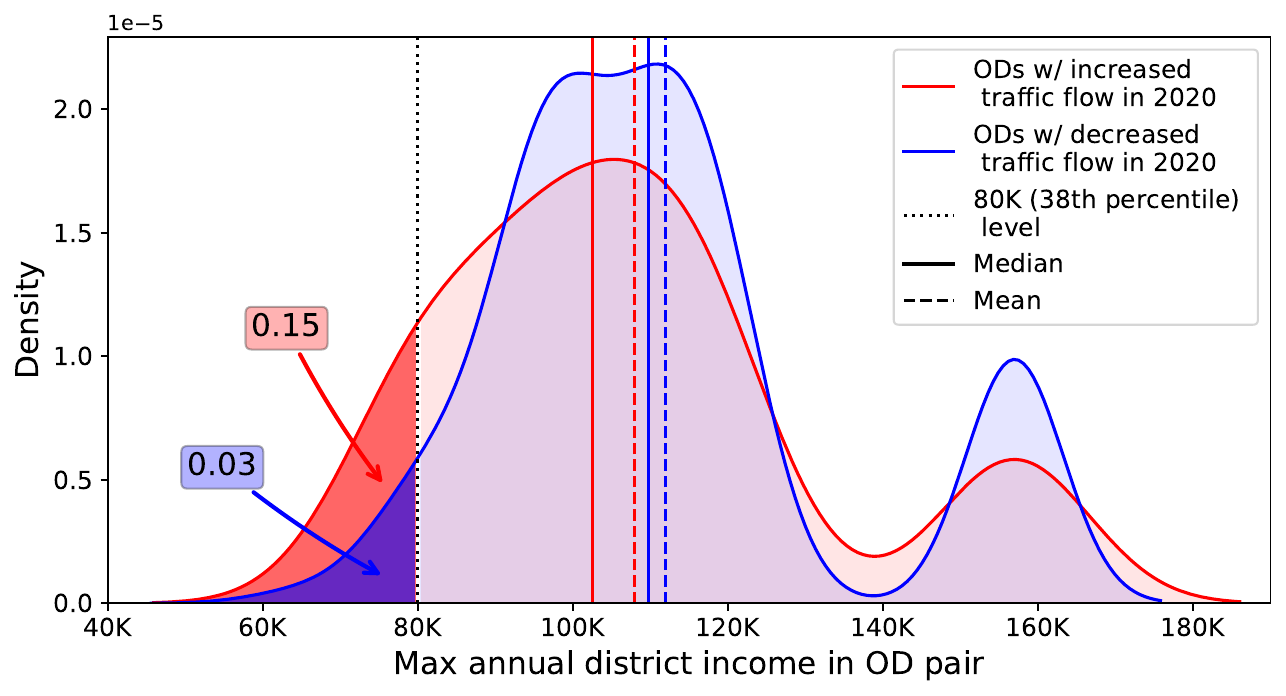}}\hfill
	\subfloat[Density of OD pairs by \emph{minimum} income in each OD pair. The dotted line represents the income at \$100,000, which corresponds to the 65th percentile of the income distribution in the health district.\label{sfig:income_min_kde}]{
		\includegraphics[width=0.7\linewidth]{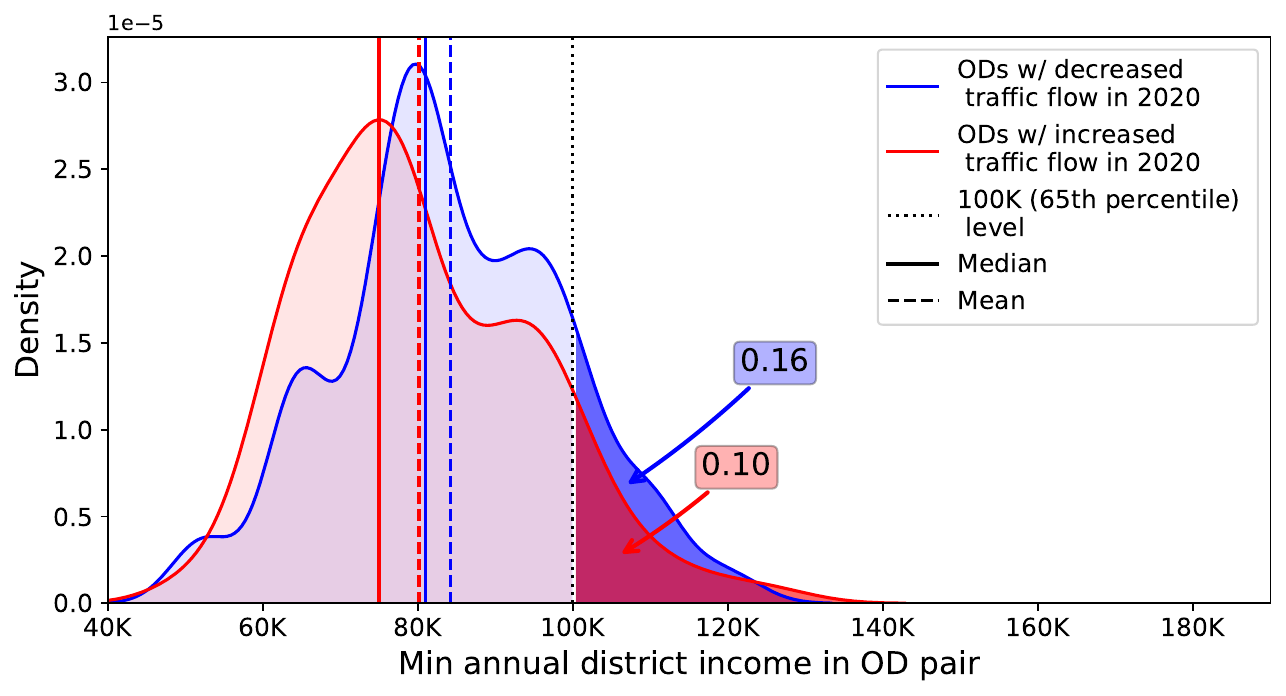}
    }\hfill
	\caption{Kernel density estimation of traffic flow changes in 2019 and 2020. The x-axis is the maximum or minimum income in each OD pair, and the y-axis is the probability density of OD pairs for each income level. The red distribution represents the kernel density across all OD pairs that showed increased traffic flow in 2020 compared to 2019. The blue distribution represents the kernel density across all OD pairs with decreased traffic flow. Both panels suggest that higher income OD pairs are more likely to have decreased traffic flow. The solid color region indicates the probability covered by this area in the distribution plot. The accompanying number represents the probability value associated with that area.}\label{fig:income_kde}
\end{figure}

In both panels of Figure \ref{fig:income_pert}, this pattern is particularly prominent in OD pairs with either an origin or destination with higher incomes (above \$90,000). This suggests that OD pairs with at least one district with higher than average income responded more to the stay-at-home order. If we analyze each individual panel in Figure \ref{fig:income_pert}, we observe distinct patterns. As shown in Figure \ref{sfig:income_max_scatter}, when the maximum annual district income in an OD pair is low (indicating travel between two low income areas), there is a higher occurrence of OD pairs with increased traffic flow. On the other hand, In Figure \ref{sfig:income_min_scatter}, when the minimum annual district income in an OD pair is high (representing travel between two high income areas), a greater number of OD pairs exhibit decreased traffic flow. This suggests that the traffic flow between low income areas tended to increase, but the traffic flow between high income areas tended to decrease during the lockdown. One possible explanation for this phenomenon is that essential workers still need to commute for work. To further analyze this trend, we examine the kernel density estimation plot illustrated in Figure \ref{fig:income_kde}.

Figure \ref{fig:income_kde} shows the kernel density estimation of OD flow by maximum or minimum income in each OD pair. In general, OD pairs with higher income are more likely to have decreased traffic flow, while those with lower income are more like to have increased flow. The x-axis is the min/max income in each OD pair, and the y-axis is the probability density of OD pairs at each income level.  The red distribution represents the kernel density across all OD pairs that showed increased traffic flow in 2020 compared to 2019. The blue distribution represents the kernel density across all OD pairs with decreased traffic flow. Both panels show that OD pairs with higher maximum or minimum income are more likely to have decreased traffic flow instead of increased flow. This is reiterated by the fact that the mean, median, and mode incomes with decreased traffic flow are higher than those with increased flow, regardless of whether we analyze this using the maximum or minimum income in each OD pair. 

Based on the findings presented in Figure \ref{fig:income_pert}, we delve deeper into the relationship between income levels and the likelihood of traffic flow increasing or decreasing. In Figure \ref{fig:income_kde}, we observe that when the maximum annual district income is low, indicating a OD traffic flow between two low income areas, the data shows a big difference between increased and decreased traffic flow. Our analysis suggests that low income areas tended to experience an increase in traffic flow during the year 2020. This conclusion is based on the findings depicted in Figure \ref{sfig:income_max_kde}, which illustrate a higher percentage of traffic flow increase compared to traffic flow decrease in regions with lower maximum annual district incomes. Specifically, when examining a threshold of \$80,000, our analysis revealed that 15\% of ODs exhibited an increase in traffic flow below this threshold. Only 3\% of traffic flow in these areas experienced a decrease. However, in Figure \ref{sfig:income_min_kde} reveals ODs are more likely to have a decreased traffic flow when the minimum annual district income is high, indicating an OD traffic flow between two high income areas. Notably, the right side of the figure displays a significant difference. When we set the threshold at \$100,000, we observe that 16\% of ODs experienced a decrease in traffic flow above this threshold, whereas 10\% exhibited an increase. On the other hand, the low income values in both panels in Figure \ref{fig:income_kde} are more likely to have increased OD flow, supporting the hypothesis that lower income employees may have been required to work more (as essential workers) or have increased travel requirements due to job-seeking as unemployment rates rose. These graphs show that OD outcomes can be useful for teasing apart the different types of effects that different households may have experienced due to COVID-19 and the stay-at-home order -- traffic to and from lower income districts may have increased, while it may have decrease in higher income districts.

\section{Conclusions}\label{sec:conclusions}

This paper presents a two stage approach to estimate 24/7 dynamic OD traffic demand using high-granular traffic counts and speed data. Specifically, we proposed incorporating additional terms into the objective function of prior models to capture urban traffic behavior, such as transitioning from local roads to highways and commuter behavior resulting in a daily ebb and flow between origins and destinations. To solve the optimization problem, we formulated the model as a non-negative least squares form and use stochastic projected gradient descent to find the optimal solution. We conducted several numerical experiments to identify a reasonable network for LA County and appropriate user weights in the optimization model between the different objective terms. The proposed data-driven optimization framework offers a promising approach for a more nuanced approach to estimating OD traffic demand in large-scale urban areas, which can be useful in a variety of applications, such as transportation planning and traffic management.

Using this formulation, we conduct a case study on traffic patterns in LA County in 2019 (reference year) and 2020 (stay-at-home order due to COVID-19). We estimated dynamic OD demands across districts in LA using high-granularity road sensor data. As expected, we found that traffic volume generally decreased during the stay-at-home order compared to 2019. However, our OD estimates provided additional nuance to these findings by allowing analysis of traffic decrease by origin-destination, which we examined by average household income in that district. We found that OD pairs with an origin or destination in a relatively higher income district decreased their traffic volume more than otherwise, suggesting that the stay-at-home order may have impacted traffic flow patterns differently across demographic groups. There could be a variety of mechanisms behind this trend; for example, commuters in higher income districts may have been able to transition to work from home during the pandemic, resulting in more dramatic reductions in OD traffic. However, this trend was not clearly the case \emph{a priori}, as the stay-at-home order also shuttered many businesses, particularly in the food and service sectors, resulting in mass layoffs that may have disproportionately affected lower income districts and their commuting patterns.

This work has several limitations, which we acknowledge here. We did not have a ground truth with which to compare OD estimates, and thus our selection of network configuration and weight parameters could not be validated precisely. However, results of our estimation process were consistent with generally expected trends in LA county. The sensor data used to parameterize the model was noisy, with occasional missing or obviously incorrect data from some sensors; while we cleaned the data to remove these values, we expect that the remaining data is subject to some measurement error. We mitigated this issue by aggregating sensor outcomes across intervals, which reduced random error. We also visually validated the data by examining daily trends, which did exhibit the appropriate peaks during rush hour periods. 

Despite these limitations, we believe that this work provides novel contributions to both the transportation and policy fields. Our extended DODE formulation is generalizable to other locations and problem applications, and we hope that this work will spur future work in both OD estimation and in using these estimates to understand sociological impacts of travel policies. By understanding the characteristics of dynamic OD demands, policymakers can make more informed decisions about transportation planning and traffic management. The estimated dynamic OD demands can also be combined with infectious disease models to predict how disease can spread between districts due to travel activity, and how these trends might be impacted by policies such as those used in LA. This information can be invaluable in predicting and preventing the spread of disease, especially during a pandemic. 
\bibliographystyle{cas-model2-names}

\bibliography{references}
\appendix
\section{Transformation to Standard L2 Norm Optimization Model}
\subsection{Base DODE Model}\label{ap:basic}
The equation \eqref{eq:twostage_final} can be reformulated in a more general form, as introduced in equation \eqref{eq:loss}. In this context, we define the decision variable OD traffic flow as $\bm q$, which is a vector of $q_{rs}^k(t)$ in the model, using the equations described in \eqref{eq:qy}. On the other hand, the input data, observed link flow $\bm y$, is defined in equation \eqref{eq:qy} by element $y_a(t)$. Here are equations.
\begin{equation}
    \begin{array}{lll}
        \arraycolsep=1.4pt\def\arraystretch{2.2}
        \bm q_{rs}(t) = \left[q_{rs}^k(t)\right]_{k\in\Phi_{rs}} ;\quad\bm q(t) = \left[\bm q_{rs}(t)\right]_{rs\in\Psi};\quad
        \bm q = \left[\bm q(t)\right]_{t\in T}^\top  \\
        \bm y(t) = \left[y_a(t)\right]_{a\in L}; \quad\bm y = \left[\bm y(t)\right]_{t\in T}^\top
    \end{array}
    \label{eq:qy}
\end{equation}
The matrix $A_{b}$ can be formulated using elements $\rho_{rs}^{ka}(t,t')$ and $p^k_{rs}(t')$, as shown in equation \eqref{eq:A}.
\begin{equation}
    \arraycolsep=1.4pt\def\arraystretch{2.2}
    \begin{array}{lll}
        \bm \rho_{rs}^{a}(t,t') = \left[\rho_{rs}^{a,1}(t,t'), ..., \rho_{rs}^{a, |\Phi_{rs}|}(t,t')\right];\quad\bm \rho^{a}(t,t') = \left[\bm \rho_{1}^{a}(t,t'), ..., \rho_{|\Psi|}^{a}(t,t')\right]\\
        \bm p_{rs}(t) =\left[ p_{rs}^1(t), ..., p_{rs}^{|\Phi_rs|}(t)\right];\quad\bm p(t) = \left[\bm p_{1}(t), ..., \bm p_{|\Psi|}(t)\right]\\
        \bm v_a(t,t') = \bm\rho^{a}(t,t') \cdot \bm p(t'); \quad\bm v_a(t) = \left[\bm v_a(t,1), ..., \bm v_a(t,|T|)\right]\\
        A_{b}(t) = \begin{bmatrix}
        \bm v_1(t)\\
        \bm v_2(t)\\
        \vdots\\
        \bm v_{|L|}(t)
        \end{bmatrix};\quad
        A_{b} = \begin{bmatrix}
        \bm A_{b}(1)\\
        \bm A_{b}(2)\\
        \vdots\\
        \bm A_{b}(|T|)
        \end{bmatrix}
    \end{array}
    \label{eq:A}
\end{equation}
Through the application of these transformations, we can readily reframe equation \eqref{eq:twostage_final} as equation \eqref{eq:1stpart}. 
\begin{equation}
    \minimize_{\bm q\geq 0} \norm{\bm y-\bm A_{b}\bm q}^2
    \label{eq:1stpart}
\end{equation}

\subsection{Incorporating local road traffic information}\label{ap:local_road}
We can formulate equation \eqref{eq:regularizer_LB} as a general regularizer $\epsilon_{LB} = \norm{\bm b_{LB} - \bm A_{LB}\bm x}$ in equation \eqref{eq:2ndformulation}. These formulations can be expressed mathematically in equations \eqref{eq:LB}.
\begin{equation}
    \arraycolsep=1.4pt\def\arraystretch{2.2}
    \begin{array}{lll}
        \bar{\bm x}(t) = \left[\bar{x}_i(t)\right]_{i\in N}; \quad\bar{\bm x} = \left[\bar{\bm x}(t)\right]_{t\in T};\quad\bm x = \left[\bm q^\top, \bar{\bm x}\right]^\top\\
        \bm b_{LB}(t) = \left[\alpha LB_i(t)\right]_{i\in N};\quad\bm b_{LB} = \left[\bm b_{LB}(t)\right]_{t\in T}^\top\\
        A_{LB}^q(i,t) = \left\{
        \arraycolsep=1.4pt\def\arraystretch{1}
        \begin{array}{ll}
            1 & \text{if there exist $\{is\}\in\Psi$  on time t}\\
            0 & \text{otherwise}
        \end{array}\right.; 
        A_{LB}^{\bar{x}}(i,t) = \left\{
        \arraycolsep=1.4pt\def\arraystretch{1}
        \begin{array}{ll}
            -1 & \text{if there exist $\{is\}\in\Psi$  on time t}\\
            0 & \text{otherwise}
        \end{array}\right.\\
        A_{LB}^q = \left[ A_{LB}^q(i,t)\right]_{i\in N, t\in T};\quad
        A_{LB}^{\bar x} = \left[ A_{LB}^{\bar x}(i,t)\right]_{i\in N, t\in T};\quad A_{LB} = \left[A_{LB}^q, A_{LB}^{\bar x}\right]
    \end{array}
    \label{eq:LB}
\end{equation}
By applying this transformation, we can rewrite equation \eqref{eq:regularizer_LB} as equation \eqref{eq:LB_final1}. 
\begin{equation}
    \begin{array}{cl}
    \epsilon_{LB}^2&=\displaystyle \sum_{i\in N}\sum_{t\in T} \norm{l_{i}(t) + d_{i}(t) - \bar{x}_i(t) - \alpha LB_i(t)}^2\\
    & = \norm{\bm b_{LB} - \bm A_{LB}^q\bm q - \bm A_{LB}^{\bar{x}}\bar{\bm x}}^2\\
    & = \norm{\bm b_{LB} - \bm A_{LB}\bm x}^2
    \end{array}
    \label{eq:LB_final1}
\end{equation}
\subsection{Incorporating daily traffic pattern information: Symmetry Constraints}\label{ap:sym}
In the next step, we can formulate equation \eqref{eq:regularizer_sym} as a general regularizer $\epsilon_{s} = \norm{\bm b_{s} - \bm A_{s}\bm x}$ in equation \eqref{eq:2ndformulation}. These formulations can
be expressed mathematically in equations \eqref{eq:sym}.
\begin{equation}
    \arraycolsep=1.4pt\def\arraystretch{2.2}
    \begin{array}{lll}
        \bar{\bm x}(t) = \left[\bar{x}_i(t)\right]_{i\in N}; \quad\bar{\bm x} = \left[\bar{\bm x}(t)\right]_{t\in T};\quad\bm x = \left[\bm q^\top, \bar{\bm x}\right]^\top\\
        \bm b_{s} = \bm 0\\
        A_{s}^q(rs,t) = \left\{
        \arraycolsep=1.4pt\def\arraystretch{1}
        \begin{array}{ll}
            1 & \text{if  $rs\in\Psi$, and $r\in R_i, s\in R_j$  on time t}\\
            -1 & \text{if $sr\in\Psi$, and $r\in R_i, s\in R_j$ on time t}\\
            0 & \text{otherwise}
        \end{array}\right.; \quad
        A_{s}^{\bar{x}} = \bm 0\\
        A_{s}^q(rs,t) = \left[ A_{s}^q(rs,t)\right]_{rs\in\Psi, t\in d, d\in D};\quad A_{s} = \left[A_{s}^q, \bm 0 \right]
    \end{array}
    \label{eq:sym}
\end{equation}

By applying this transformation, we can rewrite equation \eqref{eq:regularizer_sym} as equation \eqref{eq:sym_final1}. 
\begin{equation}
    \begin{array}{cl}
    \epsilon_{s}^2&=\displaystyle \sum_{d\in D}\sum_{R_i,R_j\in R}\norm{f_{R_i, R_j} - f_{R_j, R_i}}^2\\
    & = \norm{\bm 0 - \bm A_{s}^q\bm q}^2\\
    & = \norm{\bm b_{s} - \bm A_{s}\bm x}^2
    \end{array}
    \label{eq:sym_final1}
\end{equation}
\subsection{Imposing fidelity to total traffic flow}\label{ap:total}
In the next step, we can formulate equation \eqref{eq:regularizer_total} as a general regularizer $\epsilon_{\tau} = \norm{\bm b_{\tau} - \bm A_{\tau}\bm x}$ in equation \eqref{eq:2ndformulation}. These formulations can
be expressed mathematically in equations \eqref{eq:total}.
\begin{equation}
    \arraycolsep=1.4pt\def\arraystretch{2.2}
    \begin{array}{lll}
        \bar{\bm x}(t) = \left[\bar{x}_i(t)\right]_{i\in N}; \quad\bar{\bm x} = \left[\bar{\bm x}(t)\right]_{t\in T};\quad\bm x = \left[\bm q^\top, \bar{\bm x}\right]^\top
        \bm b_{\tau} = \left[\sum_{t\in d}\sum_{rs\in\Psi}\sum_{k\in\Phi^{rs}} \hat{q}_{rs}^{k}(t)\right]_{d\in D}\\
        A_{\tau}^q(rs,t) = \left\{
        \arraycolsep=1.4pt\def\arraystretch{1}
        \begin{array}{ll}
            1 & \text{if there exist $\{rs\}\in\Psi$  on time t}\\
            0 & \text{otherwise}
        \end{array}\right.; 
        A_{\tau}^{\bar{x}} = \bm 0\\
        A_{\tau}^q = \left[ A_{\tau}^q(rs,t)\right]_{\{rs\}\in\Psi, t\in d,d\in D};\quad A_{\tau} = \left[A_{\tau}^q, \bm 0\right]
    \end{array}
    \label{eq:total}
\end{equation}

By applying this transformation, we can rewrite equation \eqref{eq:regularizer_total} as equation \eqref{eq:total_final1}. 
\begin{equation}
    \arraycolsep=1.4pt\def\arraystretch{2.2}
    \begin{array}{cl}
    \epsilon_{\tau}^2&=\displaystyle \sum_{d\in D}\norm{\sum_{t\in d}\sum_{rs\in\Psi}\sum_{k\in\Phi^{rs}} \hat{q}_{rs}^{k}(t) - \sum_{t\in d}\sum_{rs\in\Psi}\sum_{k\in\Phi^{rs}}q^k_{rs}(t)}^2\\
    & = \norm{\bm b_{\tau} - \bm A_{\tau}^q\bm q}^2\\
    & = \norm{\bm b_{\tau} - \bm A_{\tau}\bm x}^2
    \end{array}
    \label{eq:total_final1}
\end{equation}
\end{document}